\documentclass[11pt]{article}
\usepackage{amsfonts}
\usepackage{epsfig}
\usepackage{color}
\def\be{\begin{equation}}
\def\ee{\end{equation}}
\def\bea{\begin{eqnarray}}
\def\eea{\end{eqnarray}}
\def\bes{\begin{eqnarray*}}
\def\ees{\end{eqnarray*}}
\def\pmatrix{\left(\begin{array}{cc}}
\def\endpmatrix{\end{array}\right)}

\def\sgn{{\rm sgn}}

\def\Gr{{\rm Gr}}
\def\Lag {{\rm Lag}}
\def\Gr{{\rm Gr}}

\def\det{{\rm det}}

\def\exp {{\rm exp}}

\def\Sp{{\rm Sp}}

\def\rank{{\rm rank}}

\def\nn{\nonumber}
\def\<{\langle}
\def\>{\rangle}
\def\lb{\label}
\def\bs{\setminus}

\def\R{{\bf R}}
\def\C{{\bf C}}
\def\Z{{\bf Z}}
\def\N{{\bf N}}
\def\U{{\bf U}}
\def\Q{{\bf Q}}

\def\ga{{\gamma}}

\def\om{{\omega}}
\def\Om{{\Omega}}

\def\var{{\varepsilon}}
\def\lm{{\lambda}}
\def\Lm{{\Lambda}}

\def\sg{{\sigma}}
\def\Sg{{\Sigma}}

\def\P{{\cal P}}
\def\J{{\cal J}}

\def\I{{\cal I}}
\def\E{{\cal E}}

\def\td#1{\tilde{#1}}
\def\hb{\vrule height0.18cm width0.14cm $\,$}

\def\td#1{\tilde{#1}}
\def\diag{{\rm diag}}

\title{ Seifert conjecture in the even convex case}

\author{Chungen Liu
\thanks{Partially supported by the NSF of China (11071127, 10621101), 973 Program of
MOST (2011CB808002). E-mail: liucg@nankai.edu.cn}\qquad and \qquad
Duanzhi Zhang
\thanks{Partially supported by the NSF of China
(10801078, 11171341) and Nankai University. E-mail:
zhangdz@nankai.edu.cn}\\ \\
School of Mathematics and LPMC, Nankai University\\
Tianjin 300071, People's Republic of China}

\date{}
\begin{document}
\maketitle  \begin{abstract} In this paper, we prove that there
exist at least $n$ geometrically distinct brake orbits on every
$C^2$ compact convex symmetric hypersurface $\Sg$ in $\R^{2n}$
satisfying the reversible condition $N\Sg=\Sg$ with $N=\diag
(-I_n,I_n)$. As a consequence, we show that if the Hamiltonian
function is convex and even, then Seifert conjecture of 1948 on the
multiplicity of brake orbits holds for any positive integer $n$.
\end {abstract}
\noindent {\bf Key
words:} { Brake orbit, Seifert conjecture, Lagrangian boundary
condition, Convex symmetric,   Maslov-type index}

\renewcommand{\theequation}{\thesection.\arabic{equation}}

\setcounter{equation}{0}
\section{Introduction }

For the standard symplectic space $({\bf R}^{2n},\omega_0)$ with $\omega_0(x,y)=\<Jx,y\>$, where $J=\left(\begin{array}{cc}0&-I\\I&0\\\end{array}\right)$ is
the standard symplectic matrix and $I$ is the $n\times n$ identity
matrix, an involution matrix defined by $N=\left(\begin{array}{cc}-I&0\\0&I\end{array}\right)$ is clearly anti-symplectic, i.e., $NJ=-JN$. The fixed point set of $N$ and $-N$ are the Lagrangian subspaces $L_0=\{0\}\times {\bf R}^n$  and  $L_1= {\bf R}^n\times\{0\}$ of $({\bf R}^{2n},\omega_0)$ respectively.

Suppose $H\in C^2(\R^{2n}\bs\{0\},\R)\cap C^1(\R^{2n},\R)$
satisfying   the following reversible condition
 \bea H(Nx)=H(x),\qquad
\forall\, x\in\R^{2n}.\lb{1.5}\eea We consider the following fixed
energy problem of nonlinear Hamiltonian system with Lagrangian
boundary conditions
 \bea
\dot{x}(t) &=& JH'(x(t)), \lb{1.6}\\
H(x(t)) &=& h,   \lb{1.7}\\
 x(0)&\in & L_0,\;x(\tau/2)\in L_0. \lb{1.9} \eea

\noindent It is clear that a solution $(\tau,x)$ of
(\ref{1.6})-(\ref{1.9}) is a characteristic chord on the contact
submanifold $\Sigma:=H^{-1}(h)=\{y\in\R^{2n}\,|\, H(y)=h\}$ of
$({\bf R}^{2n},\omega_0)$ and satisfies
\bea x(-t) &=& Nx(t),  \lb{1.8}\\
 x(\tau+t) &=& x(t). \lb{perio}\eea
In this paper this  kind of  $\tau$-periodic characteristic
$(\tau,x)$ is called a  {\it brake orbit} on the hypersurface $\Sg$.
We denote by ${\mathcal{J}}_b(\Sg,H)$ the set of all brake orbits on
$\Sg$. Two brake orbits $(\tau_i, x_i)\in
{\mathcal{J}}_b(\Sg,H),\;i=1,2$ are equivalent
 if the two brake orbits are
geometrically the same, i.e., $x_1(\R)=x_2(\R)$. We denote by
$[(\tau,x)]$ the equivalent class of $(\tau,x)\in
{\mathcal{J}}_b(\Sg,H)$ in this equivalent relation and by
$\td{\mathcal{J}}_b(\Sg)$ the set of $[(\tau,x)]$ for all
$(\tau,x)\in {\mathcal{J}}_b(\Sg,H)$($\td{\mathcal{J}}_b(\Sg)$ is in
fact the set of geometrically distinct brake orbits on $\Sigma$).
From now on, in the notation $[(\tau,x)]$ we always assume $x$ has
minimal period $\tau$. We also denote by $\tilde {\mathcal
{J}}(\Sg)$ the set of all geometrically distinct closed
characteristics on $\Sg$.  The number of elements in a set $S$ is
denoted by $^{\#}S$. It is well known that
$^{\#}\td{\mathcal{J}}_b(\Sg)$ (and also $^{\#}\tilde {\mathcal
{J}}(\Sg)$) is only depending on $\Sigma$, that is to say, for
simplicity we take $h=1$, if $H$ and $G$ are two $C^2$ functions
satisfying  (\ref{1.5}) and
$\Sigma_H:=H^{-1}(1)=\Sigma_G:=G^{-1}(1)$, then
$^{\#}{\mathcal{J}}_b(\Sg_H)=^{\#}{\mathcal{J}}_b(\Sg_G)$. So we can
consider the brake orbit problem in a more general setting. Let
$\Sg$ be a $C^2$ compact hypersurface in $\R^{2n}$ bounding a
compact set $C$ with nonempty interior. Suppose $\Sg$ has
non-vanishing Guassian curvature  and satisfies the reversible
condition $N(\Sg-x_0)=\Sg-x_0:=\{x-x_0|x\in \Sg\}$ for some $x_0\in
C$. Without loss of generality, we may assume $x_0=0$. We denote the
set of all such hypersurfaces  in $\R^{2n}$ by
$\mathcal{H}_b(2n)$. For $x\in \Sg$, let $N_\Sg(x)$ be the unit
outward
 normal  vector at $x\in \Sg$. Note that here by the reversible
condition there holds $N_\Sg(Nx)=NN_\Sg(x)$. We consider the
dynamics problem of finding $\tau>0$ and an absolutely continuous
curve $x:[0,\tau]\to \R^{2n}$ such that
  \bea  \dot{x}(t)&=&JN_\Sg(x(t)), \qquad x(t)\in \Sg,\lb{1.11}\\
      x(-t) &=& Nx(t), \qquad x(\tau+t) = x(t),\qquad {\rm
  for\;\; all}\;\; t\in \R.\lb{1.12}\eea

A solution $(\tau,x)$ of the problem (\ref{1.11})-(\ref{1.12})
determines  a  brake orbit on $\Sg$.

\noindent{\bf Definition 1.1.} {\it We denote  by
$$\begin{array}{ll}\mathcal{H}_b^{c}(2n)=\{\Sg\in \mathcal{H}_{b}(2n)|\;\Sg\; {\rm is\;
strictly\; convex\;} \},\\\mathcal{H}_b^{s,c}(2n)=\{\Sg\in
\mathcal{H}_{b}^c(2n)| \;-\Sg=\Sg\}.\end{array}$$}

The main result of this paper is the following

\noindent{\bf Theorem 1.1.} {\it For any
$\Sg\in\mathcal{H}_b^{s,c}(2n)$,  there holds
           $$^{\#}\td{\J}_b(\Sg)\ge n.$$}

\noindent{\bf Remark 1.1.} Theorem 1.1 is a kind of multiplicity
result related to the Arnold chord conjecture. The Arnold chord
conjecture is an existence result which was prove by K. Mohnke in
\cite{Moh}. Another kind of multiplicity result related to the
Arnold chord conjecture was proved in \cite{GL1}.

\subsection{Seifert conjecture}
Let us recall the famous conjecture proposed by H. Seifert in his
pioneer work \cite{Se1} concerning the multiplicity of brake orbits
in certain Hamiltonian systems in $\R^{2n}$.

 As a special case  of (\ref{1.5}), we assume $H\in C^2(\R^{2n},\R)$ possesses the following form
  \bea H(p,q)=\frac{1}{2}A(q)p\cdot p+V(q),\lb{i1}\eea
 where $p,q\in\R^n$, $A(q)$ is a positive definite $n\times n$ for
 any $q\in\R^n$ and $A$ is $C^2$, $V\in C^2(\R^n,\R)$ is the potential energy.
It is clear that a solution of the  following Hamiltonian system
 \bea &&\dot{x}=JH'(x),\quad x=(p,q),\lb{i2}\\&&p(0) =
p(\frac{\tau}{2})=0.\lb{i3}\eea is a brake orbit. Moreover, if $h$
is the total  energy of a brake orbit $(q, p)$, i.e.,
$H(p(t),q(t))=h$  and  $V (q(0)) = V (q(\tau )) = h$.   Then
$q(t)\in\bar\Om\equiv \{q\in \R^n|V(q)\le h\}$ for all $t\in {\bf
R}$.

In \cite{Se1} of 1948, H. Seifert studied the existence of brake
orbit for system (\ref{i2})-(\ref{i3}) with the Hamiltonian function
$H$ in the form of (\ref{i1}) and proved that
${\mathcal{J}}_b(\Sg)\neq \emptyset$ provided $V'\neq 0$ on
$\partial \Om$, $V$ is analytic and $\bar\Om$ is bounded and
homeomorphic to the unit  ball $B^n_1(0)$ in $\R^n$. Then
in the same paper he proposed the following conjecture
which is still open  for $n\ge 2$ now:\\

\qquad\qquad\qquad \qquad
$^{\#}\td{\mathcal{J}}_b(\Sg)\geq n$ {\it under the same conditions.}\\

It is well known that the lower bound $n$ in the Seifert conjecture
cannot be improved. A typical example is
 the Hamiltonian function
  \bea  H(p, q) = \frac{1}{ 2} |p|^2 + \sum_{j=1}^n a_j^2q_j^2,\qquad  q, p\in \R^n,\nn\eea
where $a_i/a_j\notin \Q$ for all $i\neq j$ and
$q=(q_1,q_2,...,q_n)$. There are exactly $n$   geometrically
distinct  brake orbits on the energy   hypersurface
$\Sigma=H^{-1}(h)$.

\subsection{Some related results since 1948}

As a special case,   letting  $A(q)=I$ in (\ref{i1}), the problem
corresponds to the following  classical fixed energy problem of the
second order autonomous Hamiltonian system
 \bea && \ddot{q}(t)+V'(q(t))=0, \quad {\rm for}\;q(t)\in\Om, \lb{1.1}\\
&& \frac{1}{2}|\dot{q}(t)|^2+V(q(t))= h, \qquad\forall t\in\R, \lb{1.2}\\
&& \dot{q}(0)=\dot{q}(\frac{\tau}{2})=0, \lb{1.3} \eea
 where $V\in C^2(\R^n, \R)$ and $h$ is  constant such that $\Om\equiv
\{q\in \R^n|V(q)<h\}$ is nonempty, bounded and connected.

 A solution $(\tau,q)$ of (\ref{1.1})-(\ref{1.3}) is still called a {\it
brake orbit$\;$} in $\bar{\Om}$. Two brake orbits $q_1$ and
$q_2:\R\to\R^n$ are {\it geometrically distinct}  if $q_1(\R)\neq
q_2(\R)$. We denote by $\mathcal{O}(\Om,V)$ and
$\td{\mathcal{O}}(\Om)$ the sets of all brake orbits and
geometrically distinct brake orbits in $\bar{\Om}$ respectively.

\noindent{\bf Remark 1.2.} It is well known that via \bea H(p,q) =
{1\over 2}|p|^2 + V(q), \nn\eea  $x=(p,q)$ and $p=\dot q$, the
elements in $ \mathcal{O}(\Omega,V)$ and the solutions of
(\ref{1.6})-(\ref{1.9}) are one to one correspondent.

\noindent{\bf Definition 1.2.}  {\it For $\Sg\in\mathcal{H}_b^{s,c}(2n)$,
a brake orbit $(\tau,x)$ on $\Sg$
 is called symmetric if
$x(\R)=-x(\R)$. Similarly, for a $C^2$ convex symmetric bounded
domain $\Omega\subset \R^n$, a brake orbit $(\tau,q)\in \mathcal
{O}(\Omega,V)$ is called symmetric if $q(\R)=-q(\R)$.}

Note that a brake orbit $(\tau,x)\in \mathcal {J}_b(\Sg,H)$ with
minimal period $\tau$
 is symmetric if
$x(t+\tau/2)=-x(t)$ for $t\in \R$, a brake orbit $(\tau,q)\in
\mathcal {O}(\Omega,V)$ with minimal period $\tau$ is  symmetric if
$q(t+\tau/2)=-q(t)$ for $t\in\R$.

After 1948, many studies have been carried out for the brake orbit
problem.
  In 1978,  S. Bolotin proved in \cite{Bol}
  the existence of brake orbits in general setting.
 In 1983-1984, K. Hayashi in \cite {Ha1}, H. Gluck and W. Ziller
 in \cite{GZ1}, and V. Benci in \cite {Be1}
  proved $^{\#}\td{\mathcal{O}}(\Om)\geq 1$ if $V$ is
$C^1$, $\bar{\Om}=\{V\leq h\}$ is compact, and $V'(q)\neq 0$ for all
$q\in \partial{\Om}$. In 1987, P. Rabinowitz in \cite{Ra1} proved
that if $H$ satisfies (\ref{1.5}), $\Sg\equiv H^{-1}(h)$ is
star-shaped, and $x\cdot H'(x)\neq 0$ for all $x\in \Sg$, then
$^{\#}\td{\mathcal{J}}_b(\Sg)\geq 1$. In 1987, V. Benci and F.
Giannoni gave a different proof of the existence of one brake orbit
in \cite{BG}. In 2005, it has been pointed out in \cite{GGP1} that
the problem of finding brake orbits is equivalent to find orthogonal
geodesic chords on manifold with concave boundary. In 2010, R.
Giamb${\rm \grave{o}}$, F. Giannoni and P. Piccione in \cite{GGP2}
proved the existence of an orthogonal geodesic chord on a Riemannian
manifold homeomorphic to a closed disk and with concave boundary.
For multiplicity of the brake problems, in 1973, A. Weinstein in
\cite{We1} proved a localized result: {\it Assume $H$ satisfies
(\ref{1.5}). For any $h$ sufficiently close to $H(z_0)$ with $z_0$
is a nondegenerate local minimum of $H$, there are $n$ geometrically
distinct brake orbits on the energy surface $H^{-1}(h)$.} In
\cite{BolZ} of 1978 and in \cite{GZ1} of 1983, under assumptions of
Seifert in \cite{Se1}, it  was proved the existence of at least $n$
brake orbits  while  a very strong assumption on the energy
integral   was used  to ensure that different minimax critical
levels correspond to geometrically distinct brake orbits.
 In 1989, A. Szulkin in \cite{Sz}
proved that $^{\#}\td{\J_b}(H^{-1}(h))\geq n$, if $H$ satisfies
conditions in \cite{Ra1} of Rabinowitz and the energy hypersurface
$H^{-1}(h)$ is $\sqrt{2}$-pinched. In 1985 E. van Groesen in
\cite{Gro} and in 1993 A. Ambrosetti, V. Benci, Y. Long in
\cite{ABL1} also proved $^{\#}\td{\mathcal{O}}(\Om)\geq n$ under
different pinching conditions. In 2006,  without pinching condition,
in \cite{LZZ} Y. Long, C. Zhu and the second author of this paper
proved that: {\it For any $\Sg\in\mathcal{H}_b^{s,c}(2n)$ with $n\ge
2$, there holds $^{\#}\td{\J}_b(\Sg)\ge 2$.}  In 2009,   the authors
of this paper in \cite{LiuZhang} proved that $^{\#}\td{\J}_b(\Sg)\ge
\left[\frac{n}{2}\right]+1$ for $\Sg\in\mathcal{H}_b^{s,c}(2n)$.
Moreover it was proved that if all brake orbits on $\Sg$ are
nondegenerate, then
         $^{\#}\td{\J}_b(\Sg)\ge n+\mathfrak{A}({\Sg}),$
where $2\mathfrak{A}(\Sigma)$ is the number of  geometrically
distinct asymmetric brake orbits on $\Sg$. Recently, in
\cite{brake1} the authors of this paper improved  the results of
\cite{LiuZhang} to that $^{\#}\td{\J}_b(\Sg)\ge
\left[\frac{n+1}{2}\right]+1$ for  $\Sg\in\mathcal{H}_b^{s,c}(2n)$,
$n\ge 3$. In \cite{brake2} the authors of this paper proved that
$^{\#}\td{\J}_b(\Sg)\ge \left[\frac{n+1}{2}\right]+2$ for
$\Sg\in\mathcal{H}_b^{s,c}(2n)$, $n\ge 4$.

\subsection{Some consequences of Theorem 1.1 and further arguments}

As  direct consequences of Theorem 1.1 we have the following two
important Corollaries.

 \noindent{\bf Corollary 1.1.} {\it If $H(p,q)$ defined by
(\ref{i1}) is even and convex, then Seifert conjecture holds.}

\noindent{\bf Remark 1.3.} If the function $H$ in Remark 1.1 is
convex and even, then $V$ is convex and even, and $\Om$ is convex
and central symmetric. Hence $\Om$ is homeomorphic to the
unit open ball in $\R^n$.

 Recently,   R. Giamb${\rm \grave{o}}$,
F. Giannoni, and P. Piccione in \cite{GGP3} gave some counterexamples to the Seifert conjecture by constructing some analytic functions $H$ with the form (\ref{i1}) such that the domain $\Omega=V^{-1}(-\infty,h)$ is homeomorphic to the unit open ball, where  $h$  is a regular value of  $V$, and there is only one brake orbit on $H^{-1}(h)$. We note that in their examples the
 functions $H$ are neither even nor convex, so  we suspect that the
convex and symmetric conditions are essential to guarantee the Seifert conjecture
in some sense.

\noindent{\bf Corollary 1.2.} {\it Suppose $V(0)=0$, $V(q)\geq 0$,
$V(-q)=V(q)$ and $V''(q)$ is positive definite for all $q\in
\R^n\bs\{0\}$. Then for any given $h>0$ and $\Om\equiv
\{q\in\R^n|V(q)<h\}$,   there holds
$$^{\#}\td{\mathcal{O}}(\Om)\ge n.$$}

It is interesting to ask the following question: {\it whether all
closed characteristics on any hypersurfaces
$\Sg\in\mathcal{H}_b^{s,c}(2n)$
 are symmetric
 brake orbits after suitable time translation provided that $^\#\td{\mathcal{J}}(\Sg)<+\infty$}?  In this direction, we have the
following result.

\noindent{\bf Theorem 1.2.} {\it For any
$\Sg\in\mathcal{H}_b^{s,c}(2n)$, suppose
         $$^{\#}\td{\J}(\Sg)= n.$$
Then all of the  $n$ closed characteristics on  $\Sigma$ are
symmetric brake orbits after suitable  time translation.}

 For $n=2$, it was proved in \cite{HWZ} that $^{\#}\tilde {\mathcal {J}}(\Sg)$ is
either $2$ or $+\infty$ for any $C^2$ compact convex hypersurface
$\Sg$ in $\R^{4}$. So Theorem 1.2 give a positive answer to the above question in the case $n=2$.
We note also that   for  the hypersurface $\Sg=\{(x_1, x_2,y_1,y_2)\in {\bf R}^4|\;x_1^2+y_1^2+\frac{x_2^2+y_2^2}{4}=1\}$ there hold $^\#\td{\mathcal{J}_b}(\Sg)=+\infty$ and $^\#\td{\mathcal{J}^s_b}(\Sg)=2$. Here we denote  by $\td{\mathcal{J}^s_b}(\Sg)$ the set of all symmetric brake orbits on $\Sg$. We also note that on the hypersurface $\Sg=\{x\in {\bf R}^{2n}|\;|x|=1\}$ there are some non-brake closed characteristics.

The key ingredients in the proof of Theorem 1.1 are some ideas from
our previous paper \cite{LiuZhang} and the following result
which generalizes corresponding results of our previous papers
\cite{brake1,brake2} completely, where the iteration path $\gamma^2$
will be defined in Definition 2.5 below.

\noindent{\bf Theorem 1.3.} {\it For $\ga\in \P_\tau(2n)$, let
$P=\ga(\tau)$. If $i_{L_0}(\gamma)\ge 0$, $i_{L_1}(\gamma)\ge 0$, $i(\ga)\ge n$,
$\ga^2(t)=\ga(t-\tau)\ga(\tau)$ for all $t\in[\tau,2\tau]$,  then
 \be
 i_{L_1}(\ga)+S_{P^2}^+(1)-\nu_{L_0}(\ga)\ge 0.\lb{pp}\ee}

In this paper, we denote by $\N$, $\Z$, $\Q$, $\R$ and $\C$ the sets
of positive integers, integers, rational numbers, real numbers and
complex numbers respectively. We denote by both $\langle
\cdot,\cdot\rangle$ and $\cdot$ the standard inner product in $\R^n$
or $\R^{2n}$, by $(\cdot,\cdot)$ the inner product of corresponding
Hilbert space. For any $a\in \R$, we denote by $[a]=\sup\{k\in
\Z|k\le a\}$.

\setcounter{equation}{0}\section{   Index theories for symplectic paths and the homotopic properties of  symplectic matrices }   

In this section we make some preparations for the proof of Theorems
1.1-1.2. We first briefly introduce the Maslov-type index theory of
$(i_{L_j},\nu_{L_j})$ for $j=0,1$ and $(i_\om,\nu_\om)$ for
$\om\in\U:=\{z\in\C|\,|z|=1\}$.

Let $ \mathcal {L}(\R^{2n})$ denotes the set of $2n\times 2n$  real
matrices  and $\mathcal{L}_s(\R^{2n})$ denotes its subset of
symmetric ones. For any $F\in \mathcal{L}_s(\R^{2n})$, we denote by
$m^*(F)$ the dimension of maximal positive
     definite subspace, negative definite subspace, and kernel of
     any $F$ for $*=+,-,0$ respectively.

Let $J_k=\left(\begin{array}{cc}0&-I_k\\I_k&0\end{array}\right)$ and
$N_k=\left(\begin{array}{cc}-I_k&0\\0&I_k\end{array}\right)$ with
$I_k$ being the identity in $\R^k$. If $k=n$ we will omit the
subscript $k$ for convenience, i.e., $J_n=J$ and $N_n=N$.

 The symplectic group $\Sp(2k)$ for any $k\in\N$ is
defined by
 $$\Sp(2k)=\{M\in \mathcal{L}(\R^{2k})| M^TJ_kM=J_k\},$$
where $M^T$ is the transpose of matrix $M$.

For any $\tau>0$, the symplectic path in $\Sp(2k)$ starting from the
identity $I_{2k}$ is defined by
  $$\mathcal{P}_{\tau}(2k)=\{\gamma\in C([0,\tau],\Sp(2k))|
  \gamma(0)=I_{2k}\}.$$

 In the study of periodic solutions of Hamiltonian
systems, the Maslov-type index theory of $(i(\gamma),\nu(\gamma))$
of $\gamma$ usually plays a important role which was introduced by
C. Conley and E. Zehnder in \cite{CoZ} for nondegenerate symplectic
path $\gamma\in\mathcal{ P}_\tau(2n)$ with $n\ge 2$, by Y. Long and
E. Zehnder in \cite{LZe} for nondegenerate symplectic path
$\gamma\in\mathcal {P}_\tau(2)$, by Long in \cite{Long4} and C.
Viterbo in \cite{V} for $\gamma\in\mathcal {P}(2n)$. In
\cite{Long0}, Long introduced the $\omega$-index which is an index
function $(i_{\omega}(\gamma),\nu_{\omega}(\gamma))\in \Z\times
\{0,1,\cdots,2n\}$ for $\omega\in \U$.

 For any
$\om\in \U$, the following hypersurface in $\Sp(2n)$ is defined by:
    $$ \Sp(2n)_\om^0=\{M\in \Sp(2n)|\det(M-\om I_{2n})=0\}.$$
For any two continuous path $\xi$ and $\eta$: $[0,\tau]\to\Sp(2n)$
with $\xi(\tau)=\eta(0)$, their joint path is defined by
    \be  \eta * \xi(t)=\left\{\begin{array}{lr} \xi(2t)\qquad &{\rm
     if}\,
     0\le t\le \frac{\tau}{2},\\ \eta(2t-\tau)\quad &{\rm if}\,
     \frac{\tau}{2}\le t\le \tau .\end{array}\right.\lb{zzz}\ee
Given any two $(2m_k\times 2m_k)$  matrices of square block form
$M_k=\left(\begin{array}{cc}A_k&B_k\\C_k&D_k\end{array}\right)$ for
$k=1,2$, as in \cite{Long1}, the $\diamond$-product (or symplectic direct product) of $M_1$ and
$M_2$ is defined by the following $(2(m_1+m_2)\times 2(m_1+m_2))$
matrix $M_1\diamond M_2$:
     $$M_1\diamond M_2=\left(\begin{array}{cccc}A_1&0&B_1&0\\0&A_2&0&B_2\\
        C_1&0&D_1&0\\0&C_2&0&D_2\end{array}\right).$$
 We denote by $M^{\diamond k}$ the $k$-times self
$\diamond$-product of $M$ for any $k\in \N$.

 It is easy to see that
 $$N_{m_1+m_2}(M_1\diamond M_2)^{-1}N_{m_1+m_2}(M_1\diamond M_2)
   =(N_{m_1}M_1^{-1}N_{m_1}M_1)\diamond
   (N_{m_2}M_2^{-1}N_{m_2}M_2).$$
A special path $\xi_n$ is defined by
        $$\xi_n(t)=\left(\begin{array}{cc}2-\frac{t}{\tau}&0\\0&(2-\frac{t}{\tau})^{-1}\end{array}\right)^{\diamond n}, \qquad \forall t\in[0,\tau].$$
{\bf Definition 2.1.} {\it For any $\om\in\U$ and $M\in\Sp(2n)$, define
     \be \nu_\om(M)=\dim_\C\ker(M-\om I_{2n}).\ee
For any $\ga\in \mathcal{P}_\tau(2n)$, define
  \be \nu_\om(\ga)=\nu_\om(\ga(\tau)).\ee
If $\ga(\tau)\notin\Sp(2n)_\om^0$, we define
     \be i_\om(\ga)=[\Sp(2n)_\om^0\,:\,\ga*\xi_n],\lb{m1}\ee
where the right-hand side of (\ref{m1}) is the usual homotopy
intersection number and the orientation of $\ga*\xi_n$ is its
positive time direction under homotopy with fixed endpoints. when
$\om=1$ we will write $i_1(\ga)$ as $i(\ga)$ in convenience. If
$\ga(\tau)\in\Sp(2n)_\om^0$, we let $\mathcal{F}(\ga)$ be the set of
all open neighborhoods of $\ga$ in $\mathcal{P}_\tau(2n)$, and
define
 \be i_\om(\ga)=\sup_{U\in\mathcal{F}(\ga)}\inf\{i_\om(\beta)|\,\beta(\tau)\in U \,{\rm and}\,
 \beta(\tau)\notin\Sp(2n)_\om^0\}.\ee}
 The index pair $(i_\om(\ga),\nu_\om(\ga))\in \Z\times \{0,1,...,2n\}$, which is
called the index function of $\ga$ at $\om$,  was first defined in a
different way by Y. Long in \cite{Long0}(see also \cite{Long1} and
\cite{Long4}).

 For any $M\in\Sp(2n)$ we define
 \bea \Om(M)=\{P\in \Sp(2n)&|&\sg(P)\cap\U=\sg(M)\cap \U\nn\\
                            && {\rm and}\,
                            \nu_\lm(P)=\nu_\lm(M),\;\; \forall
                            \lm\in\sg(M)\cap\U\},\eea
where we denote by $\sg(P)$ the spectrum of $P$.

 We denote by $\Om^0(M)$ the path connected component of
$\Om(M)$ containing $M$, and call it the {\it homotopy component} of
$M$ in $\Sp(2n)$.

\noindent{\bf Definition 2.2.} {\it For any $M_1$,$M_2\in \Sp(2n)$, we
call $M_1\approx M_2$ if $M_1\in \Om^0(M_2)$.}

\noindent{\bf Remark 2.1.} It is easy to check that $\approx$ is an
equivalent relation. If $M_1\approx M_2$, we have $M_1^k\approx
M_2^k$ for any $k\in \N$ and $M_1\diamond M_3 \approx M_2\diamond
M_4$ for $M_3\approx M_4$. Also we have $M_1\diamond M_2\approx
M_2\diamond M_1$ and $PMP^{-1}\approx M$ for any $P,M\in \Sp(2n)$.
By Theorem 7.8 of \cite{Long0}, $M_1\diamond M_2\approx M_1 \diamond
M_3$ if and only if $M_2\approx M_3$.

\noindent{\bf Lemma 2.1.} {\it Assume $M_1\in \Sp(2(k_1+k_2))$
and $M_2\in\Sp(2k_3)$ have the following block form
$M_1=\left(\begin{array}{cccc}A_1&A_2&B_1&B_2\\A_3&A_4&B_3&B_4\\C_1&C_2&D_1&D_2\\C_3&C_4&D_3&D_4\end{array}\right)$
and $M_2=\left(\begin{array}{cc} A_5&B_5\\C_5&D_5\end{array}\right)$
with $A_1, B_1,C_1,D_1\in \mathcal{L}(\R^{k_1})$, $A_4,
B_4,C_4,D_4\in \mathcal{L}(\R^{k_1})$, $A_5,D_5\in
\mathcal{L}(\R^{k_3})$. Let
$M_3=\left(\begin{array}{cccccc}A_1&0&A_2&B_1&0&B_2\\0&A_5&0&0&B_5&0\\
     A_3&0&A_4&B_3&0&B_4\\C_1&0&C_2&D_1&0&D_2\\0&C_5&0&0&D_5&0\\
     C_3&0&C_4&D_3&0&D_4\end{array}\right)$. Then
     \bea M_3\approx M_1\diamond M_2.\lb{pp1}\eea}

{\bf Proof.} Let
$P=\diag\left(\left(\begin{array}{ccc}I_{k_1}&0&0\\0&0&I_{k_2}\\0&I_{k_3}&0\end{array}\right),
\left(\begin{array}{ccc}I_{k_1}&0&0\\0&0&I_{k_2}\\0&I_{k_3}&0\end{array}\right)\right)$.
It is east to verify that $P\in\Sp(2(k_1+k_2+k_3))$ and
$M_3=P(M_1\diamond M_2)P^{-1}$. Then (\ref{pp1}) holds from Remark
2.1 and the proof of Lemma 2.1 is complete.\hfill\hb

The following symplectic matrices were introduced as {\it basic
normal forms} in \cite{Long1}:
      \bea D(\lm)=
      \left(\begin{array}{cc}\lm&0\\0&\lm^{-1}\end{array}\right),\qquad
     && \lm=\pm 2,\nn\\
      N_1(\lm,b)=\left(\begin{array}{cc}\lm&b\\0&\lm\end{array}\right),\qquad
     && \lm=\pm1,\,b=\pm1,\,0,\nn\\
      R(\theta)=\left(\begin{array}{cc}\cos(\theta)&-\sin(\theta)\\
      \sin(\theta)&\cos(\theta)\end{array} \right), \qquad
      &&\theta\in(0,\pi)\cup(\pi,2\pi),\nn\\
      N_2(\om,b)=\left(\begin{array}{cc}R(\theta)&b\\0&R(\theta)\end{array}\right),
        \qquad
      &&\theta\in(0,\pi)\cup(\pi,2\pi),\nn\eea
where $b= \left(\begin{array}{cc}b_1&b_2\\b_3&b_4\end{array}\right)$
with $b_i\in\R$ and $b_2\neq b_3$.

For any $M\in \Sp(2n)$ and $\om\in\U$, {\it splitting number} of $M$
at $\om$ is defined by
   \bea S_M^{\pm}(\om)=\lim_{\epsilon\to 0^+}
   i_{\om\exp(\pm\sqrt{-1}\epsilon)}(\ga)-i_\om(\ga)\nn\eea
for any path $\ga\in\mathcal{P}_\tau(2n)$ satisfying $\ga(\tau)=M$.

Splitting numbers possesses the following properties.

 \noindent{\bf Lemma 2.2.} ( \cite{Long0}, Lemma 9.1.5 and List
9.1.12 of \cite{Long1}) {\it Splitting number $S_M^{\pm}(\om)$ are
well defined, i.e., they are independent of the choice of the path
$\ga\in\mathcal{P}_\tau(2n)$ satisfying $\ga(\tau)=M$. For
$\om\in\U$ and $M\in\Sp(2n)$, $S_Q^{\pm}(\om)=S_M^{\pm}(\om)$ if
$Q\approx M$. Moreover we have

(1) $(S_M^+(\pm1),S_M^-(\pm1))=(1,1)$ for $M=\pm N_1(1, b)$ with
$b=1$ or $0$;

(2) $(S_M^+(\pm1),S_M^-(\pm1))=(0,0)$ for $M=\pm N_1(1, b)$ with
$b=-1$;

 (3) $(S_M^+(e^{\sqrt{-1}\theta}),S_M^-(e^{\sqrt{-1}\theta}))=(0,1)$ for
   $M=R(\theta)$ with $\theta\in(0,\pi)\cup (\pi,2\pi)$;

(4) $(S_M^+(\om),S_M^-(\om))=(0,0)$
 for $\om\in \U\setminus \R$ and $M=N_2(\om, b)$ is {\bf trivial} i.e., for

  sufficiently small
 $\alpha>0$, $MR((t-1)\alpha)^{\diamond n}$ possesses no eigenvalues
 on $\U$ for $t\in [0,1)$.

(5) $(S_M^+(\om),S_M^-(\om)=(1,1)$
 for $\om\in \U\setminus \R$ and $M=N_2(\om, b)$ is {\bf non-trivial}.

 (6) $(S_M^+(\om),S_M^-(\om)=(0,0)$
for any $\om\in\U$ and $M\in \Sp(2n)$ with $\sg(M)\cap
\U=\emptyset$.

 (7) $S_{M_1\diamond M_2}^\pm(\om)=S_{M_1}^\pm(\om)+S_{M_2}^\pm(\om)$, for any $M_j\in \Sp(2n_j)$ with
    $j=1,2$ and $\om\in\U$.}

Let \be F=\R^{2n}\oplus \R^{2n}\lb{zhang0}\ee possess the standard
inner product. We define the symplectic structure of $F$ by
      \be \{v,w\}=(\mathcal{J}v,w),\;\forall v,w\in F,\;
 {\rm where}\; \mathcal{J}=(-J)\oplus J=\left(\begin{array}{cc} -J &0\\0&J\end{array}\right).
 \;\lb{zhang2}\ee
  We denote by $\Lag(F)$ the set of Lagrangian subspaces of $F$,
  and equip it with the topology as a subspace of the Grassmannian of all
  $ 2n$-dimensional subspaces of $F$.

 It is easy to check that, for any $M\in \Sp(2n)$ its
  graph
      $$\Gr(M)\equiv\left\{\left(\begin{array}{c}x\\Mx\end{array}\right)|x\in
      \R^{2n}\right\}$$
is a Lagrangian subspace of $F$.

Let \bea V_1= L_0\times L_0=\{0\}\times \R^n\times \{0\}\times
\R^n\subset \R^{4n},\\ \quad V_2= L_1\times L_1=\R^n\times
\{0\}\times \R^n\times\{0\}\subset \R^{4n}.\lb{s2}\eea

By Proposition 6.1 of \cite{LZ} and Lemma 2.8 and Definition 2.5 of
\cite{LZZ}, we give the following

 \noindent{\bf Definition 2.3.} {\it For any continuous path $\ga\in\mathcal{P}_\tau(2n)$, we
 define the following Maslov-type indices:

 \bea && i_{L_0}(\ga)=\mu^{CLM}_{F}(V_1, \Gr(\ga),[0,\tau])-n,\\
     && i_{L_1}(\ga)=\mu^{CLM}_{F}(V_2, \Gr(\ga),[0,\tau])-n,\\
      && \nu_{L_j}(\ga)=\dim (\ga(\tau)L_j\cap L_j),\qquad j=0,1,\eea
where we denote by $i^{CLM}_F(V,W,[a,b])$ the Maslov index for
Lagrangian subspace path pair $(V,W)$ in $F$ on $[a,b]$ defined by
Cappell, Lee, and Miller in \cite{CLM}. For any $M\in\Sp(2n)$ and
$j=0,1$, we also denote by $\nu_{L_j}(M)=\dim(ML_j\cap L_j)$.}

The index $i_{L}(\ga)$ for any Lagrangian subspace $L\subset {\bf
R}^{2n}$ and symplectic path $\ga\in \mathcal P_{\tau}(2n)$ was
defined by the first author of this paper in \cite{Liu2} in a
different way(see also \cite{Liu0} and \cite{LZZ}).

\noindent{\bf Definition 2.4.}  {\it For two paths
$\gamma_0,\;\gamma_1\in
 \mathcal{P_\tau}(2n)$ and $j=0,1$, we say that they are $L_j$-homotopic and denoted by
 $\gamma_0\sim_{L_j}\gamma_1$, if there is a map $\delta:[0,1]\to
 \mathcal{P}(2n)$ such that $\delta(0)=\gamma_0$ and $\delta(1)=\gamma_1$, and
 $\nu_{L_j}(\delta(s))$ is constant for $s\in [0,1]$.}

  \noindent{\bf Lemma 2.3.}(\cite{Liu2}) {\it (1) If  $\gamma_0\sim_{L_j}\gamma_1$, there hold
  $$i_{L_j}(\gamma_0)=i_{L_j}(\gamma_1),\;\nu_{L_j}(\gamma_0)=\nu_{L_j}(\gamma_1).$$
 (2) If $\gamma=\gamma_1\diamond \gamma_2\in
  \mathcal{P}(2n)$, and correspondingly $L_j=L_j'\oplus L_j''$, then
  $$i_{L_j}(\gamma)=i_{L'_j}(\gamma_1)+i_{L_j''}(\gamma_2),\;\nu_{L_j}(\gamma)=\nu_{L'_j}(\gamma_1)+\nu_{L_j''}(\gamma_2).$$
 (3) If
$\gamma\in \mathcal{P}(2n)$ is the fundamental solution of
$$\dot x(t)=JB(t)x(t)$$ with symmetric matrix function
$B(t)=\pmatrix b_{11}(t) & b_{12}(t)\\b_{21}(t) & b_{22}(t)
\endpmatrix$ satisfying $b_{22}(t)>0$ for any $t\in R$, then there holds
$$i_{L_0}(\gamma)=\sum_{0<s<1}\nu_{L_0}(\gamma_s),\;\gamma_s(t)=\gamma(st).$$
(4) If $b_{11}(t)>0$ for any $t\in \R$, there holds
$$i_{L_1}(\gamma)=\sum_{0<s<1}\nu_{L_1}(\gamma_s),\;\gamma_s(t)=\gamma(st).$$  }

\noindent {\bf Definition 2.5.} {\it For any $\ga\in
\mathcal{P}_\tau$ and $k\in \N\equiv\{1,2,...\}$, in this paper the
$k$-time iteration $\ga^k$ of $\ga\in \mathcal{P}_\tau(2n)$ in brake
orbit boundary sense is defined by $\td{\ga}|_{[0,k\tau]}$ with \bea
   \td{\ga}(t)=\left\{\begin{array}{l}\ga(t-2j\tau)(N\ga(\tau)^{-1}N\ga(\tau))^j,
   \; t\in[2j\tau,(2j+1)\tau], j=0,1,2,...\\
      N\ga(2j\tau+2\tau-t)N(N\ga(\tau)^{-1}N\ga(\tau))^{j+1}\; t\in[(2j+1)\tau,(2j+2)\tau],
      j=0,1,2,...\end{array}\right.\nn\eea}

\setcounter{equation}{0}
\section {$(L_0,L_1)$-concavity and $(\varepsilon, L_0,L_1)$-signature of symplectic matrix}  

\noindent {\bf Definition 3.1.} {\it For any $P\in\Sp(2n)$ and
$\varepsilon\in\R$,   we define the $(\var,L_0,L_1)$-symmetrization
of $P$ by
 \bea M_\varepsilon(P)=P^T\left(\begin{array}{cc}\sin{2\var}I_n&-\cos{2\var I_n}\\-\cos{2\var}I_n&-\sin
 2\var I_n
 \end{array}\right)P+\left(\begin{array}{cc}\sin{2\var}I_n&\cos{2\var}I_n\\\cos{2\var}I_n&-\sin2\var
 I_n
 \end{array}\right).\nn\eea
  The $(\varepsilon, L_0,L_1)$-signature of $P$ is defined by the
signature of $M_\varepsilon(P)$. The $(L_0,L_1)$-concavity and
$(L_0,L_1)^*$-concavity of a symplectic path $\gamma$ is defined by
$$concav_{(L_0,L_1)}(\ga)=i_{L_0}(\ga)-i_{L_1}(\ga),\;\;
concav_{(L_0,L_1)}^*(\ga)=(i_{L_0}(\ga)+\nu_{L_0}(\ga))
-(i_{L_1}(\ga)+\nu_{L_1}(\ga))$$ respectively.}

 In \cite{Liu2}  it
was proved that $(L_0,L_1)$-concavity is only depending on the
 end matrix $\gamma(\tau)$ of $\gamma$,  and in \cite{Zhang1} it was proved that the $(L_0,L_1)$-concavity of a symplectic path
  $\gamma$ is a half of the $(\varepsilon, L_0,L_1)$-signature of $\gamma(\tau)$. i.e., we have the following result.

\noindent{\bf Theorem 3.1.}{ (\cite{Zhang1})}  {\it For
$\ga\in\mathcal{P}_\tau(2k)$ with $\tau>0$, we have
   \bea concav_{(L_0,L_1)}(\ga)=\frac{1}{2}\sgn
M_\var(\ga(\tau)),\nn\eea where $\sgn M_\var(\ga(\tau))$ is the
signature of the symmetric matrix $M_\var(\ga(\tau))$ and $0<\var\ll
1$. we also have,
 \bea
 concav_{(L_0,L_1)}^*(\ga)=\frac{1}{2} \sgn
M_\var(\ga(\tau)),\;\;0<-\var\ll 1.\nn\eea}

\noindent{\bf Remark 3.1.} (Remark 2.1 of \cite{Zhang1}) {\it For
any $n_j \times n_j$ symplectic matrix $P_j$ with $j=1,2$ and
 $n_j\in\N$, we have
    \bea &&M_\var(P_1\diamond P_2)=M_\var(P_1)\diamond
    M_\var(P_2),\nn\\
     &&\sgn M_\var(P_1\diamond P_2)=\sgn M_\var(P_1)+
    \sgn M_\var(P_2),\nn\eea
  where $\var\in\R$.}

 In the rest of this section, we further develope some basic
properties of the $(\varepsilon, L_0,L_1)$-signature
 and study the normal form of $L_0$-degenerate symplectic matrices.

\noindent{\bf Lemma 3.1.} (Lemma 2.3 of \cite{brake1}) {\it Let
$k\in\N$ and any symplectic matrix
$P=\left(\begin{array}{cc}I_k&0\\C&I_k\end{array}\right)$. Then
$P\approx I_2^{\diamond p}\diamond N_1(1,1)^{\diamond q}\diamond
N_1(1,-1)^{\diamond r}$ with $p=m^0(C),\;q=m^-(C),\;r=m^+(C)$.}

\noindent{\bf Definition 3.2.} {\it We call two symplectic matrices
  $M_1$ and $M_2$ are $(L_0,L_1)$-{\it homotopic equivalent} in $\Sp(2k)$, and denote it by
$M_1\sim M_2$, if there are $P_j\in\Sp(2k)$ with
$P_j=\diag (Q_j,(Q_j^T)^{-1})$, where $Q_j$ is a $k\times k$
invertible real matrix, and $\det(Q_j)>0$ for $j=1,2$, such that
    $$M_1=P_1M_2P_2.$$}
 {\bf Remark 3.1.} Let $M_i=\left(\begin{array}{cc}A_i & B_i\\C_i & D_i\end{array}\right)\in \Sp(2k_i)$, $i=0,1,2$ and $M_1\sim M_2$($k_1=k_2$ in this time), then
$A_1^TC_1$, $B_1^TD_1$ are congruent to $A_2^TC_2$, $B_2^TD_2$
respectively. So $m^*(A_1^TC_1)=m^*(A_2^TC_2)$ and
$m^*(B_1^TD_1)=m^*(B_2^TD_2)$ for $*=\pm,\;0$. Furthermore, if
$M_0=M_1\diamond M_2$(here $k_1=k_2$ is not necessary), then \bea
m^*(A_0^TC_0)=m^*(A_1^TC_1)+m^*(A_2^TC_2),\quad
m^*(B_0^TD_0)=m^*(B_1^TD_1)+m^*(B_2^TD_2).\lb{gogo}\eea
So $m^*(A^TC)$ and $m^*(B^TD)$ are $(L_0,L_1)$-homotopic  invariant. The
following formula will be used frequently
 \bea N_kM_1^{-1}N_kM_1=I_{2k}+2\left(\begin{array}{cc}B_1^TC_1 & B_1^TD_1\\A_1^TC_1 & C_1^TB_1\end{array}\right).\lb{use}\eea

It is clear that $\sim$ is an equivalent relation and we have the
following

\noindent{\bf Lemma 3.2.} (Lemma 2.4 of \cite{brake1}) {\it For
$M_1,\,M_2\in \Sp(2k)$, if $M_1\sim M_2$, then
 \bea
        \sgn M_\var(M_1)=\sgn M_\var(M_2),\quad 0\le |\var|\ll 1,\nn\\
       N_kM_1^{-1}N_kM_1\approx N_kM_2^{-1}N_kM_2.\nn\eea}
By results in \cite{brake1, brake2, Zhang1}, we have the following
lemmas 3.3-3.5 which will be used frequently in Section 4.

 \noindent{\bf Lemma 3.3.} ( Lemma 2.5 of \cite{brake1}) {\it
Assume $P=\left(\begin{array}{cc}A&B\\C&D\end{array}\right)\in
\Sp(2k)$, where $A,B,C,D$ are all $k\times k$ matrices.

(i)  Let  $q=\max\{m^+(A^TC),m^+(B^TD)\}$, we have
      \bea \frac{1}{2}\sgn M_\var(P)\le k-q-\nu_{L_1}(P),\quad 0<-\var\ll 1,\nn\eea
     \bea \frac{1}{2}\sgn M_\var(P)\le k-q-\nu_{L_0}(P),\quad 0<\var\ll 1.\nn\eea

(ii) If both $B$ and $C$ are invertible, we have
       \bea \sgn M_\var(P)=\sgn M_0(P),\quad 0\le |\var|\ll 1.\nn\eea}

\noindent{\bf Lemma 3.4.} (\cite{Zhang1}) {\it For
$\ga\in\mathcal{P}_\tau(2)$, $b>0$, and $\var>0$ small enough we
have
 \bea
&& \sgn M_{\pm\var}(R(\theta))=0,\quad
 {\rm for}\;\theta\in \R,\nn\\
 &&  \sgn M_{\pm\var}(P)=0,\quad
 {\rm if}\;P=\left(\begin{array}{cc} a & 0\\0 & 1/a\end{array}\right) \;{\rm with}\;a\in \R\setminus\{0\},\nn\\
&&  \sgn M_\var(P)=0,\quad {\rm if}\; P=\pm
 \left(\begin{array}{cc}1&b\\0&1\end{array}\right)\;{\rm or}\;
 \pm
 \left(\begin{array}{cc}1&0\\-b&1\end{array}\right),\nn\\
&&\sgn M_\var(P)=2,\quad {\rm if}\; P=\pm
 \left(\begin{array}{cc}1&-b\\0&1\end{array}\right),\nn\\
 &&\sgn M_\var(P)=-2,\quad {\rm if}\; P=
 \pm
 \left(\begin{array}{cc}1&0\\b&1\end{array}\right).\nn\eea}\\
\noindent{\bf Lemma 3.5.} (Lemma 2.9 of \cite{brake2}) {\it Let $ 2
k\times 2 k$ symmetric real matrix $E$ have the following block form
$E =\left(\begin{array}{cc}0&E_1\\E_1^T&E_2\end{array}\right)$. Then
     \be m^\pm(E)\ge \rank E_1.\ee}\\
 In the following we prove Lemma 3.6, which will be used to prove
Lemma 3.7 while Lemma 3.7 and Lemma 3.8 are two key lemmas in this
paper.

 \noindent{\bf Lemma 3.6.} {\it Let $A_1$ and $A_3$ be $k\times
k$ real matrices. Assume both $A_1$ and $A_1A_3$ are symmetric and
$\sg(A_3)\subset (-\infty,0)$. Then we have
       \bea \sgn A_1+\sgn (A_1A_3)=0.\lb{n1}\eea}
\qquad{\bf Proof.} It is clear that $A_3$ is invertible. We prove
Lemma 3.6 by the following two steps.

\noindent{\bf Step 1.} We prove this lemma in the case $A_1$ is
invertible by mathematical induction for $k\in \N$.

If $k=1$, then $A_1, A_3\in \R$ and (\ref{n1}) holds obviously. Now
assume (\ref{n1}) holds for $1\le k\le l$. If we can prove
(\ref{n1}) for $k=l+1$, then by the mathematical induction
(\ref{n1}) holds for any $k\in\N$ and Lemma 3.6 is proved in the
case $A_1$ is invertible.

By the real Jordan canonical form decomposition of $A_3$, in Step 1
we only need to prove (\ref{n1}) for $k=l+1$ in the following Case 1
and Case 2.

\noindent{\bf Case 1.} There is an invertible $(l+1)\times (l+1)$
real
 matrix such that $Q^{-1}A_3Q$ is the $(l+1)$-order Jordan form
  $\left(\begin{array}{cccccc} \lm&1&0&\cdots &0&0\\0&\lm&1&\cdots&0&0\\
    \vdots&\vdots&\vdots&\cdots&\vdots&0\\0&0&0&\cdots&\lm&1\\0&0&0&\cdots&0&\lm
    \end{array}\right):=\td{A}_3$ with $\lm<0$.

Denote by $\td {A}_1=Q^TA_1Q$. We have
 \bea \td{A}_1\td{A}_3=Q^TA_1Q\,Q^{-1} A_3Q=Q^TA_1A_3Q.\nn\eea
So both $\td{A}_1$ and $\td{A}_1\td{A}_3$ are symmetric and we have
 \bea \sgn A_1+\sgn (A_1A_3)=\sgn \td{A}_1+\sgn
 (\td{A}_1\td{A}_3).\lb{n2}\eea
Denote by $\td{A}_1=(a_{i,j})$, where $a_{i,j}$ is the element on
the $i$-th row and $j$-th column of $\td{A}_1$ for $1\le i,j\le
l+1$. We denote by $\td{A}_1\td{A}_3=(c_{i,j})$ in the same sense.
Then we have $a_{i,j}=a_{j.i}$ and $c_{i,j}=c_{j,i}$ for $1\le
i,j\le l+1$.

\noindent{\bf Claim 3.1.} In Case 1 $a_{i,j}=0$ for $i+j\le l+1$ and
$a_{i,j}=a_{l+1,1}$ for $i+j=l+2$ with $1\le i,j\le l+1$.

 For $2\le j\le l+1$, since $c_{1,j}=c_{j,1}$ we have
  \bea \lm a_{1,j}+a_{1,j-1}=\lm a_{j,1}=\lm a_{1,j}.\nn\eea
  So we have
  \bea a_{1,j-1}=0,\quad  2\le j\le l+1.\lb{n3}\eea
For $2\le i,j\le l+1$, since $c_{i,j}=c_{j,i}$ we have
 \bea \lm a_{i,j}+a_{i,j-1}=\lm a_{j,i}+a_{j,i-1}=\lm a_{i,j}+a_{i-1,j}.\nn\eea
  So we have
  \bea a_{i,j-1}=a_{i-1,j},\quad 2\le i, j\le l+1.\lb{n4}\eea
By (\ref{n3}) and (\ref{n4}) we have
  \bea
 && a_{i,j}=a_{i-1,j+1}=\cdots=a_{2,i+j-2}=a_{1,i+j-1}=0,
  \quad 1\le i,j \;{\rm and}\; i+j\le l+1,\lb{n5}\\
  && a_{l+1,1}=a_{l,2}=a_{l-1,3}=\cdots=a_{2,l}=a_{1,l+1}.\lb{n6}
  \eea
Then Claim 3.1 holds from (\ref{n5}) and (\ref{n6}).

By Claim 3.1, let $a=a_{1,l+1}$ we have
\bea \td{A}_1=\left(\begin{array}{ccccccc} 0&0&0&0&0&0&a\\0&0&0&0&0&a&*\\
    0&0&0&0&\cdot&*&*\\0&0&0&\cdot&*&*&*\\0&0&\cdot&*&*&*&*\\0&a&*&*&*&*&*
    \\a&*&*&*&*&*&*
    \end{array}\right),\quad \td{A}_1\td{A}_3=\left(\begin{array}{ccccccc} 0&0&0&0&0&0&\lm a\\0&0&0&0&0&\lm a&*\\
    0&0&0&0&\cdot&*&*\\0&0&0&\cdot&*&*&*\\0&0&\cdot&*&*&*&*\\0&\lm a&*&*&*&*&*
    \\\lm a&*&*&*&*&*&*
    \end{array}\right).\eea
Then it is easy to see that $\td{A}_1\td{A}_3$ is congruent to $\lm
\td{A}_1$. So since $\lm<0$ we have
 \bea\sgn
 (\td{A}_1\td{A}_3)=\sgn(\lm\td{A}_1)=-\sgn(\td{A}_1).\nn\eea
 Hence we have
  \bea \sgn
 (\td{A}_1\td{A}_3)+\sgn\td{A}_1=0.\lb{n7}\eea
Then (\ref{n1}) holds from (\ref{n2}) and (\ref{n7}). So in Case 1
(\ref{n1}) holds for $k=l+1$.

\noindent{\bf Case 2.} There exists a invertible $(l+1)\times (l+1)$
real matrix $Q$ such that $Q^{-1}A_3Q=\diag(A_4,A_5)$, where $A_4$
is a $k_1\times k_1$ real matrix with $\sg({A_4})\subset
(-\infty,0)$ and $A_5$ is a $k_2$-order Jordan form
$$ A_5=\left(\begin{array}{cccccc} \lm&1&0&\cdots &0&0\\0&\lm&1&\cdots&0&0\\
    \vdots&\vdots&\vdots&\cdots&\vdots&0\\0&0&0&\cdots&\lm&1\\0&0&0&\cdots&0&\lm
    \end{array}\right)$$ with $\lm<0$, $1\le k_1,k_2\le l$ and $k_1+k_2=l+1$.

We still denote by $\td {A}_1=Q^TA_1Q$, then
 \bea \td{A}_1\td{A}_3=Q^TA_1Q\,Q^{-1} A_3Q=Q^TA_1A_3Q.\nn\eea
So both $\td{A}_1$ and $\td{A}_1\td{A}_3$ are symmetric and we have
 \bea \sgn A_1+\sgn (A_1A_3)=\sgn \td{A}_1+\sgn
 (\td{A}_1\td{A}_3).\lb{n8}\eea
Correspondingly we can write $\td{A}_1$ in the block form
decomposition
$\td{A}_1=\left(\begin{array}{cc}E_1&E_2\\E_2^T&E_4\end{array}\right)$,
where $E_1$ is a $k_1\times k_1$ real symmetric matric and $E_4$ is
a $k_2\times k_2$ real symmetric matrix. Then
 \bea
 \td{A}_1\td{A}_3=\left(\begin{array}{cc}E_1A_4&E_2A_5\\E_2^TA_4&E_4A_5\end{array}\right)\nn\eea
is symmetric.

\noindent {\bf Subcase 1.} $E_4$ is invertible.

In this case we have
   \bea && \left(\begin{array}{cc}I_{k_1}&-E_2E_4^{-1}\\0&I_{k_2}\end{array}\right)
         \left(\begin{array}{cc}E_1&E_2\\E_2^T&E_4\end{array}\right)
         \left(\begin{array}{cc}I_{k_1}&0\\-E_4^{-1}E_2^T&I_{k_2}\end{array}\right)\nn\\
        &=&\left(\begin{array}{cc}E_1-E_2E_4^{-1}E_2^T&0\\0&E_4\end{array}\right)\lb{n9}\eea
and
     \bea &&\left(\begin{array}{cc}I_{k_1}&-E_2E_4^{-1}\\0&I_{k_2}\end{array}\right)
         \left(\begin{array}{cc}E_1A_4&E_2A_5\\E_2^TA_4&E_4A_5\end{array}\right)
         \left(\begin{array}{cc}I_{k_1}&0\\-E_4^{-1}E_2^T&I_{k_2}\end{array}\right)\nn\\
        &=&\left(\begin{array}{cc}E_1A_4-E_2E_4^{-1}E_2^TA_4&0\\0&E_4A_5\end{array}\right)\nn\\
        &=&\left(\begin{array}{cc}(E_1-E_2E_4^{-1}E_2^T)A_4&0\\0&E_4A_5\end{array}\right).\lb{n10}\eea
Since $\td{A}_1$ is symmetric and invertible, by (\ref{n9}) we have
$E_1-E_2E_4^{-1}E_2^T$ is symmetric and invertible. Since
$\td{A}_1\td{A}_3$ is symmetric and invertible, by {\ref{n10}) we
have $(E_1-E_2E_4^{-1}E_2^T)A_4$ is symmetric and invertible. Since
$1\le k_1\le l$ and $\sg(A_4)\subset (-\infty,0)$, by our induction
hypothesis we have
 \bea \sgn ((E_1-E_2E_4^{-1}E_2^T)A_4)+\sgn
 (E_1-E_2E_4^{-1}E_2^T)=0.\lb{n11}\eea
By (\ref{n10}) we also have $E_4A_5$ is symmetric. Since $E_4$ is
symmetric and invertible, $\sg(A_5)=\{\lm\}\subset(-\infty,0)$ and
$1\le k_2\le l$, by our induction hypothesis we have
    \bea \sgn (E_4A_5)+\sgn E_4=0.\lb{n12}\eea
By (\ref{n9}) we have
 \bea \sgn \td{A}_1=\sgn (E_1-E_2E_4^{-1}E_2^T)+\sgn
 E_4.\lb{n13}\eea
By (\ref{n10}) we have
   \bea \sgn(\td{A}_1\td{A}_3)=\sgn ((E_1-E_2E_4^{-1}E_2^T)A_4)+\sgn
   (E_4A_5).\lb{n14}\eea
Then by (\ref{n11})-(\ref{n14}) we have
 \bea \sgn(\td{A}_1\td{A}_3)+\sgn \td{A}_1=0.\lb{n15}\eea
Then (\ref{n1}) holds from (\ref{n8}) and (\ref{n15}).

 \noindent{\bf
Subcase 2.} $E_4$ is not invertible.

In this case we define $k_2$-order real invertible matrix
    \bea E_0=\left(\begin{array}{ccccccc} 0&0&0&0&0&0&1\\0&0&0&0&0&1&0\\
    0&0&0&0&\cdot&0&0\\0&0&0&\cdot&0&0&0\\0&0&\cdot&0&0&0&0\\0&1&0&0&0&0&0
    \\1&0&0&0&0&0&0
    \end{array}\right).\nn\eea
Then  it is easy to verify that $E_0A_5$ is symmetric and
$E_4+\var E_0$ is invertible for $0<\var \ll 1$. Define
$A_\var=\left(\begin{array}{cc}E_1&E_2\\E_2^T&E_4+\var
E_0\end{array}\right)$. Since $\td{A}_1$ and $\td{A}_1\td{A}_3$ are
invertible, we have both $A_\var$ and $A_\var\td{A}_3$ are symmetric
and invertible. So we have
      \bea \sgn \td{A}_1=\sgn A_\var,\quad \sgn
      (\td{A}_1\td{A}_3)=\sgn(A_\var\td{A}_3), \quad {\rm for}\; 0<\var\ll 1.\lb{n16}\eea
By the proof of Subcase 1, we have
     \bea \sgn(A_\var\td{A}_3)+\sgn A_\var=0.\lb{n17}\eea
So by ({\ref{n16}) we have
    \bea \sgn(\td{A}_1\td{A}_3)+\sgn \td{A}_1=0.\lb{n18}\eea
Then (\ref{n1}) holds from (\ref{n18}).

So in Case 2 (\ref{n1}) holds for $k=l+1$. Hence in the case $A_1$
is invertible Lemma 3.6.holds and Step 1 is finished.

\noindent{\bf Step 2.} We prove (\ref{n1}) in the case $A_1$ is not
invertible.

   If $A_1=0$, (\ref{n1}) holds obviously.

   If $1\le \rank A_1=m\le k-1$, there is a real orthogonal matrix $G$ such that
       \bea G^TA_1G=\left(\begin{array}{cc}
       0&0\\0&\hat{A}_1\end{array}\right),\lb{n23}\eea
where $\hat{A}_1$ is a $m$-order invertible real symmetric matrix.
Correspondingly we write
 \bea G^{-1}A_3G=\left(\begin{array}{cc}
       F_1&F_2\\ F_3&F_4\end{array}\right),\nn\eea
 where $F_1$ is a $(k-m)\times(k-m)$ real matrix and $F_4$ is a
 $m\times m$ real matrix.

 Since $A_1A_3$ is symmetric, we have
  \bea G^TA_1A_3G=G^TA_1GG^{-1}A_3G=\left(\begin{array}{cc}
       0&0\\\hat{A}_1 F_3&\hat{A}_1F_4\end{array}\right)\nn\eea
is still symmetric. So we have $\hat{A}_1F_2^T=0$, since $\hat{A}_1$
is invertible we have $ F_3=0$. Then
    \bea G^{-1}A_3G=\left(\begin{array}{cc}
       F_1&  F_2 \\0&F_4\end{array}\right).\lb{n20}\eea
So we have
       \bea G^TA_1A_3G^T=\left(\begin{array}{cc}
       0&0\\0&\hat{A}_1F_4\end{array}\right)\lb{n21}\eea
and $\hat{A}_1F_4$ is symmetric. Also by (\ref{n20}) we have $F_4$
is invertible and $\sg(F_4)\subset(-\infty,0)$. So by the proof of
the case $A_1$ is invertible we have
     \bea \sgn (\hat{A}_1F_4)+\sgn \hat{A}_1=0.\lb{n22}\eea
By (\ref{n23}) and (\ref{n21}) we have
  \bea \sgn(A_1A_3)+\sgn A_1=\sgn (\hat{A}_1F_4)+\sgn
  \hat{A}_1.\lb{n24}\eea
Then (\ref{n1}) holds from (\ref{n22}) and (\ref{n24}). Hence
(\ref{n1}) holds in the case $A_1$ is not invertible. Step 2 is
finished.

By Step 1 and Step 2  Lemma 3.6 holds.\hfill\hb

\noindent{\bf Lemma 3.7.} {\it Let
$R=\left(\begin{array}{cc}A_1&I_k\\A_3&A_2\end{array}\right)\in\Sp(2k)$
with $A_3$  being  invertible. If $e(N_kR^{-1}N_kR)=2m$, where
$0\le m\le k$ and the elliptic hight $e(P)$ of $P$ is the total
algebraic multiplicity of all eigenvalues of $P$ on $\bf U$ for any
$P\in\Sp(2n)$. Then we have
      \bea m-k\le \frac{1}{2}\sgn M_\var(R)\le k-m,\qquad 0\le |\var|\ll
       1.\lb{n25}\eea}

{\bf Proof.} Since $e(N_kR^{-1}N_kR)=2m$, there exists a symplectic
matrix $P\in\Sp(2k)$ such that
 \bea P^{-1}(N_kR^{-1}N_kR)P=Q_1\diamond Q_2\lb{n28}\eea
with $\sg(Q_1)\in \U$, $\sg(Q_2)\cap\U=\emptyset$, $Q_1\in\Sp(2m)$,
and $Q_2\in\Sp(2k-2m)$. By (ii) of Lemma 3.3, since $A_3$ is
invertible we only need to prove (\ref{n25}) for $\var=0$.

\noindent{\bf Step 1.} We first prove (\ref{n25}) in the case $A_1$
is invertible.

 Since $R$ is a symplectic matrix we have
   $R^TJ_kR=J_k$. Then  $A_1^TA_3$ and $A_2$ are all
   symmetric matrices and
    \bea A_1^TA_2-A_3^T=I_k.\nn\eea
 Since $R^T$ is also a symplectic matrix we have
$RJ_kR^T=J_k$. Then $A_1$ is symmetric. Hence $A_1A_3$ is symmetric
and
  \be A_1A_2-A_3^T=I_k.\lb{liu3}\ee
 By
  definition we have
      \bea
      M_0(R)&=& R^T\left(\begin{array}{cc}0&-I_k\\-I_k&0\end{array}\right)R+
       \left(\begin{array}{cc}0&I_k\\I_k&0\end{array}\right)\nn\\
            &=&-2\left(\begin{array}{cc}A_1 A_3&A_3^T\\A_3&A_2\end{array}\right).\eea
Since $A_1$ is invertible, we have
     \bea && \left(\begin{array}{cc}I_k&0\\-A_1^{-1}&I_k\end{array}\right)
     \left(\begin{array}{cc}A_1A_3&A_3^T\\A_3&A_2\end{array}\right)
     \left(\begin{array}{cc}I_k&-A_1^{-1}\\0&I_k\end{array}\right)\nn\\
      &=&\left(\begin{array}{cc}A_1A_3&0\\0&-A_1^{-1}A_3^T+A_2\end{array}\right)\nn\\
      &=&\left(\begin{array}{cc}A_1A_3&0\\0&A_1^{-1}\end{array}\right),\lb{liu4}\eea
    where in the last equality we have used  the equality (\ref{liu3}).
So by (\ref{liu4}) we have
 \bea
\frac{1}{2}\sgn M_0(R)=-\frac{1}{2}\sgn
  \left(\begin{array}{cc}
  A_1A_3&0\\0& A_1^{-1}\end{array}\right).\lb{n26}\eea

By the Jordan canonical form decomposition of complex matrix, there
exists a complex invertible $k$-order matrix $G_1$ such that
\bea G_1^{-1}A_3G_1=\left(\begin{array}{ccccc} u_1&*&*&*&*\\0&u_2&*&*&*\\
    0&0&\ddots&*&*\\0&0&0&u_{k-1}&*\\0&0&0&0&u_k
    \end{array}\right)\nn\eea
    with $u_1,u_2,...,u_k\in\C$.

By (\ref{use}) we have
 \bea
N_kR^{-1}N_kR=I_{2k}+2\left(\begin{array}{cc}A_3&A_2\\A_1A_3&A_3^T\end{array}\right).\lb{liu7}\eea
Since
    \bea \left(\begin{array}{cc}I_k&0\\-A_1&I_k\end{array}\right)
         \left(\begin{array}{cc}A_3&A_2\\A_1A_3&A_3^T\end{array}\right)
          \left(\begin{array}{cc}I_k&0\\A_1&I_k\end{array}\right)\nn
          =\left(\begin{array}{cc}I_k+2A_3&A_2\\-A_1&-I_k\end{array}\right),\eea
by (\ref{liu7}) we have
      \bea
      \left(\begin{array}{cc}I_k&0\\A_1&I_k\end{array}\right)^{-1}(N_kR^{-1}N_kR)
       \left(\begin{array}{cc}I_k&0\\A_1&I_k\end{array}\right)=
       \left(\begin{array}{cc}3I_k+4A_3&2A_2\\-2A_1&-I_k\end{array}\right):=R_1.\lb{liu8}\eea
By (\ref{liu8}), for any $\lm\in\C$ we have
       \be \lm I_{2k}-R_1=\left(\begin{array}{cc}(\lm-3)
       I_k-4A_3&-2A_2\\2A_1&(\lm+1)I_k\end{array}\right).\lb{liu9}\ee
 Since $A_1$ is invertible, by (\ref{liu3}) we have
     \bea &&\left(\begin{array}{cc}I_k&-\frac{1}{2}((\lm-3)I_k-4A_3)A_1^{-1}\\0&I_k\end{array}\right)
            \left(\begin{array}{cc}(\lm-3)
       I_k-4A_3&-2A_2\\2A_1&(\lm+1)I_k\end{array}\right)\nn\\
        &=&\left(\begin{array}{cc}0&-\frac{1}{2}((\lm^2-2\lm+1)I_k-4\lm A_3)A_1^{-1}
        \\2A_1&(\lm+1)I_k\end{array}\right).\lb{liu10}\eea
   Then by (\ref{liu9})-(\ref{liu10}) we have
    \bea \det (\lm I_{2k}-R_1)
    =\det ((\lm^2-2\lm+1)I_k-4\lm A_3).\lb{liu12}\eea
    Denote
    by $u_1,u_2,...,u_k$ the $k$ complex eigenvalues of $A_3$, by (\ref{liu12}) we have
   \bea \det (\lm I_{2k}-R_1)=\Pi_{i=1}^k (\lm^2-2\lm+1-4\lm u_i)=\Pi_{i=1}^k (\lm^2-(2+4u_i)\lm+1).\lb{liu13}\eea
So by (\ref{liu8}) and (\ref{liu13}) we have
   \bea \det (\lm I_{2k}-N_kR^{-1}N_kR)=\Pi_{i=1}^k (\lm^2-2\lm+1-4\lm u_i)=
   \Pi_{i=1}^k (\lm^2-(2+4u_i)\lm+1).\lb{n27}\eea
It is easy to check that the equation $\lm^2-(2+u_i)\lm+1=0$ has two
solutions on $\U$ if and only if $-4\le u_i\le 0$ for
$i=1,2,3...,k$. So by (\ref{n28}) without loss of generality we
assume $u_j \in [-4,0)$ for $1\le j\le m$ and $u_j\notin [-4,0)$ for
$m+1\le j\le k$. Then there exists a real invertible matrix
$k$-order Q such that
 \bea Q^{-1}A_3Q=\left(\begin{array}{cc}A_4&0\\0&A_5\end{array}\right):=\td{A}_3\nn\eea
and $\sg(A_4)\subset [-4,0)$, $\sg(A_5)\cap [-4,0)=\emptyset$, where
$A_4$ is an $m$-order real invertible matrix and $A_5$ is a
$(k-m)$-order real matrix.

Denote by $\td {A}_1=Q^TA_1Q$. We have
 \bea \td{A}_1\td{A}_3=Q^TA_1Q\,Q^{-1} A_3Q=Q^TA_1A_3Q.\nn\eea
So both $\td{A}_1$ and $\td{A}_1\td{A}_3$ are symmetric and we have
 \bea \sgn A_1+\sgn (A_1A_3)=\sgn \td{A}_1+\sgn
 (\td{A}_1\td{A}_3).\lb{n29}\eea
Correspondingly we can write $\td{A}_1$ in the block form
decomposition
$\td{A}_1=\left(\begin{array}{cc}E_1&E_2\\E_2^T&E_4\end{array}\right)$,
where $E_1$ is an $m$-order real symmetric matric and $E_4$ is a
$(k-m)$-order real symmetric matrix. Then
 \bea
 \td{A}_1\td{A}_3=\left(\begin{array}{cc}E_1A_4&E_2A_5\\E_2^TA_4&E_4A_5\end{array}\right)\nn\eea
is symmetric.

By the same argument of the proof of Subcase 2 of Lemma 3.6 without
loss of generality we can assume $E_1$ is invertible(Otherwise we
can perturb it slightly such that it is invertible). So as in
Subcase 1 of the proof of Lemma 3.6 we have
   \bea && \left(\begin{array}{cc}I_{m}&0\\-E_2^TE_1^{-1}&I_{k-m}\end{array}\right)
         \left(\begin{array}{cc}E_1&E_2\\E_2^T&E_4\end{array}\right)
         \left(\begin{array}{cc}I_{m}&-E_1^{-1}E_2\\0&I_{k-m}\end{array}\right)\nn\\
        &=&\left(\begin{array}{cc}E_1&0\\0&E_4-E_2^TE_1^{-1}E_2\end{array}\right)\lb{n30}\eea
and
     \bea &&\left(\begin{array}{cc}I_{m}&0\\-E_2^TE_1^{-1}&I_{k-m}\end{array}\right)
         \left(\begin{array}{cc}E_1A_4&E_2A_5\\E_2^TA_4&E_4A_5\end{array}\right)
         \left(\begin{array}{cc}I_{m}&-E_1^{-1}E_2\\0&I_{k-m}\end{array}\right)\nn\\
                &=&\left(\begin{array}{cc}E_1A_4&0\\0&(E_4-E_2^TE_1^{-1}E_2)A_5\end{array}\right).\lb{n31}\eea
 By (\ref{n31})
we also have $E_1A_4$ is symmetric. Since $E_1$ is symmetric and
invertible, $\sg(A_4)\subset [-4,0)$, by Lemma 3.6 we have
    \bea \sgn (E_1A_4)+\sgn E_1=0.\lb{n32}\eea
By (\ref{n30}) we have
 \bea \sgn \td{A}_1=\sgn (E_4-E_2^TE_1^{-1}E_2)+\sgn
 E_1.\lb{n33}\eea
By (\ref{n31}) we have
   \bea \sgn(\td{A}_1\td{A}_3)=\sgn ((E_4-E_2^TE_1^{-1}E_2)A_5)+\sgn
   (E_1A_4).\lb{n34}\eea
Then by (\ref{n32})-(\ref{n34}) we have
 \bea \begin{array}{ll} \sgn(\td{A}_1\td{A}_3)+\sgn \td{A}_1\\ =\sgn ((E_4-E_2^TE_1^{-1}E_2)A_5)
 +\sgn (E_4-E_2^TE_1^{-1}E_2)\in [-2(k-m),2(k-m)].\end{array}\lb{n35}\eea
Then (\ref{n25}) holds from (\ref{n26}), (\ref{n29}) and
(\ref{n35}).

\noindent{\bf Step 2.} We prove (\ref{n25}) in the case $A_1$ is not
invertible.

If $A_1=0$, then $A_3=-I_k$ and $m=k$. It is easy to check that
$M_0(R)=2\left(\begin{array}{cc}0 & I_k\\I_k &
-A_2\end{array}\right)$ is congruent to $2\left(\begin{array}{cc}0 &
I_k\\I_k & 0\end{array}\right)$, so $\sgn M_0(R)=0$, (\ref{n25})
holds.

If $1\le \rank A_1=r\le k-1$, there is a $k\times k$ invertible
matrix $G$ with $\det G>0$ such that
     \be (G^{-1})^TA_1G^{-1}=\diag(0,\Lm),\ee
where $\Lm$ is a $r\times r$ real invertible matrix. Then we have
 \bea  \diag((G^T)^{-1},G)\cdot R\cdot\diag(G^{-1},G^T)&=&
 \left(\begin{array}{cc}(G^T)^{-1}A_1G^{-1}&I_k\\GA_3G^{-1}&GA_2G^T\end{array}\right) \nn\\&:
 =&R_2=
 \left(\begin{array}{cccc}0&0&I_{k-r}&0\\0&\Lm&0&I_r
 \\
 B_1&B_2&D_1&D_2\\B_3&B_4&D_3&D_4\end{array}\right),\lb{liu14}\eea
where $B_1$ and $D_1$ are $(k-r)\times(k-r)$ matrices and $B_4$ and
$D_4$ are $r\times r$ matrices.

Then since $R_2$ is symplectic and $\Lm$ is invertible, we have
 $R_2^TJ_kR_2=J_k$ which implies that $B_3=0$, $D_3=D_2^T$,
  $B_1=-I_{k-r}$, and $D_1$, $D_4$ are symmetric. So we have
 \bea R_2=\left(\begin{array}{cccc}0&0&I_{k-r}&0\\0&\Lm&0&I_r
 \\
 B_1&B_2&D_1&D_2\\0&B_4&D_2^T&D_4\end{array}\right)\nn\eea
For $t\in[0,1]$, we define
 \bea \beta(t)=\left(\begin{array}{cccc}0&0&I_{k-r}&0\\0&\Lm&0&I_r
 \\
 B_1&tB_2&tD_1&tD_2\\0&B_4&tD_2^T&D_4\end{array}\right)\nn\eea
Then it is easy to check that $\beta$ is a symplectic path and
$\nu_{L_j}(\beta(t)=0$ for all $t\in[0,1]$ and $j=0,1$. Also we have
$\beta(1)=R_2$ and
 \bea \beta(0)=\left(\begin{array}{cccc}0&0&I_{k-r}&0\\0&\Lm&0&I_r
 \\ B_1&0&0&0\\0&B_4&0&D_4\end{array}\right)=-J_{k-r}\diamond\left(\begin{array}{cc}\Lm&I_r\\
     B_4&D_4\end{array}\right):=R_3.\nn\eea

Then by Lemma 2.2 of \cite{Zhang1}, Lemma 3.4, and Remark 3.1 we
have
  \bea \frac{1}{2}\sgn M_0(R_2)&=&\frac{1}{2}\sgn
  M_0(-J_{k-r})+\frac{1}{2}\sgn M_0 \left(\left(\begin{array}{cc}\Lm&I_r\\
     B_4&D_4\end{array}\right)\right)\nn\\&=&\frac{1}{2}\sgn M_0 \left(\left(\begin{array}{cc}\Lm&I_r\\
     B_4&D_4\end{array}\right)\right).\lb{liu30}\eea
Since $R_2\sim R$, by (\ref{liu30}) we have
 \bea \frac{1}{2}\sgn M_0(R)=\frac{1}{2}\sgn M_0 \left(\left(\begin{array}{cc}\Lm&I_r\\
     B_4&D_4\end{array}\right)\right).\lb{liu31}\eea

By (\ref{use}) we have
     \bea N_kR_2^{-1}N_kR_2=I_{2k}+2\left(\begin{array}{cccc}B_1&B_2&D_1&D_2\\0&B_4&D_2^T&D_4
 \\ 0&0&B_1^T&0\\0&\Lm B_4& B_2^T&B_4^T\end{array}\right).\lb{liu15}\eea

By (\ref{liu15}) for any $\lm\in\C$, we have
 \bea &&\det(\lm I_{2k}-N_kR_2^{-1}N_kR_2)\nn\\&=&
      \det((\lm-1) I_{k-r}-2B_1)\det((\lm-1)
      I_{k-r}-2B^T_1)\cdot\nn\\&&\cdot
       \det\left(\begin{array}{cc}(\lm-1)I_r-2B_4&-2D_4\\-2\Lm
       B_4&(\lm-1)I_r-2B_4^T\end{array}\right)\nn\\
       &=&\det(\lm I_{2k}-N_kR_3^{-1}N_kR_3),\lb{liu16}\eea
where
  \bea N_kR_3^{-1}N_kR_3=I_{2k}+2\left(\begin{array}{cccc}B_1&0&0&0\\0&B_4&0&D_4
 \\ 0&0&B_1^T&0\\0&\Lm B_4&0&B_4^T\end{array}\right).\nn\eea
So by (\ref{liu16}) we have
     \bea\sg(N_kR^{-1}N_kR)=\sg(N_kR_2^{-1}N_kR_2)=\sg(N_kR_3^{-1}N_kR_3).\lb{liu17}\eea
 Since $B_1=-I_{k-r}$ and $R_3=(-J_{k-r})\diamond\left(\begin{array}{cc}\Lm&I_r\\
     B_4&D_4\end{array}\right)$, by (\ref{liu17}) we have
  \bea e\left(N_r\left(\begin{array}{cc}\Lm&I_r\\
     B_4&D_4\end{array}\right)^{-1}N_r\left(\begin{array}{cc}\Lm&I_r\\
     B_4&D_4\end{array}\right)\right)=2(m-(k-r)).\lb{liu18}\eea
So by step 1 we have
  \bea \frac{1}{2}\left|\sgn M_0 \left(\left(\begin{array}{cc}\Lm&I_r\\
     B_4&D_4\end{array}\right)\right)\right|\le
     r-(m-(k-r))=k-m.\lb{liu19}\eea
   Then (\ref{n25}) holds form (\ref{liu31}) and (\ref{liu19}). Thus
   Step 2 is finished.

  By Step 1 and Step 2, Lemma 3.7 holds.
   \hfill\hb

The following result is about the $(L_0,L_1)$-normal form of
$L_0$-degenerate symplectic matrices which generalizes Lemma 2.10 of
\cite{brake2}.

\noindent {\bf Lemma 3.8.}   {\it Let
$R\in\Sp(2k)$ has the block form
$R=\left(\begin{array}{cc}A&B\\C&D\end{array}\right)$ with
$1\le\rank B=r<k$. We have

(i) $R\sim
\left(\begin{array}{cccc}A_1&B_1&I_r&0\\0&D_1&0&0\\A_3&B_3&A_2&0\\
C_3&D_3&C_2&D_2
\end{array}\right)$, where $A_1,A_2,A_3$ are $r\times r$ matrices,
$D_1,D_2,D_3$ are $(k-r)\times (k-r)$ matrices, $B_1,B_3$ are
$r\times (k-r)$ matrices, and $C_2,C_3$ are $(k-r)\times r$
matrices.

(ii) If $A_3$ is invertible, we have
     \be R\sim \left(\begin{array}{cc}A_1&I_r\\A_3&A_2\end{array}\right)\diamond
     \left(\begin{array}{cc}D_1&0\\\td{D}_3&D_2\end{array}\right),\ee
     where $\td{D}_3$ is a $(k-r)\times (k-r)$ matrix.

(iii) If $1\le \lb=\rank A_3=\lm\le r-1$, then
     \be R\sim
            \left(\begin{array}{cc}U&I_\lm\\\Lm &
       V\end{array}\right)\diamond\left(\begin{array}{cccc}
     \td{A}_1&\td{B}_1&I_{r-\lm}&0\\
       0&D_1&0&0\\0&\td{B}_3&\td{A}_2&0\\
       \td{C}_3&\td{D}_3&\td{C}_2&\td{D}_2\end{array}\right),\ee
     where $\td{A}_1,\td{A}_2$ are $(r-\lm)\times (r-\lm)$ matrices,
     $\td{B}_1,\td{B}_3$ are $(r-\lm)\times (k-r)$ matrices,
     $\td{C}_2, \td{C}_3$ are $(k-r)\times (r-\lm)$ matrices,
     $D_1,\td{D}_2,\td{D}_3$ are $(k-r)\times (k-r)$ matrices,
     $U, V, \Lm$ are $\lm\times\lm$ matrices, and $\Lm$ is invertible.

 (iv) If $A_3=0$,  then $A_1$, $A_2$ are symmetric and $A_1A_2=I_r$.   Suppose $m^+(A_1)=p$, $m^-(A_1)=r-p$ and  $0\le\rank B_3=\lm\le\min\{r,\,k-r\}$, then
  \bea  &&N_kR^{-1}N_kR\approx  \left(\begin{array}{cc}1&1\\0&1\end{array}\right)^{\diamond p+q^-}
  \diamond \left(\begin{array}{cc}1&-1\\0&1\end{array}\right)^{\diamond (r-p+q^+)
  }\diamond I_2^{\diamond q^0}\diamond
      D(2)^{\diamond \lm},\lb{m1}\\
      &&m^+(A^TC)=\lm+ q^+,\lb{m2}\\ &&m^0(A^TC)=r-\lm+q^0,\lb{m3}\\&& m^-(A^TC)=\lm+q^-,\lb{m4}\eea
   where $q^*\ge 0$ for $*=\pm,0$, $q^++q^0+q^-=k-r-\lm$,  for any symplectic matrix the term $M^{\diamond 0}$ means it does not appear}.

{\bf Proof.} By Lemma 2.10 of \cite{brake2} or the same argument of
the proof of Theorem 3.1 of \cite{brake1}, (i)-(iii) hold. So we
only need to prove (\ref{m1})-(\ref{m4}).

By (i) and $A_3=0$ we have
 \bea R\sim
\left(\begin{array}{cccc}A_1&B_1&I_r&0\\0&D_1&0&0\\0&B_3&A_2&0\\
C_3&D_3&C_2&D_2
\end{array}\right):=R_1.\lb{m12}\eea

\noindent Since $R_1$ is symplectic we have $R_1^TJ_kR_1=J_k$. Then
we have $A_1$, $A_2$ are symmetric and $A_1A_2=I_r$.
$D_1D_2^T=I_{k-r}$ and $A_1^TB_3=C_3^TD_1$. By (\ref{use}) we have
    \bea N_kR_1^{-1}N_kR_1=\left(\begin{array}{cccc}
       I_r&2B_3&2A_2&0\\0&I_{k-r}&0&0\\0&2A_1^TB_3&I_r&0\\2B_3^TA_1&2B_1^TB_3
       +2D_1^TD_3&2B^T_3&I_{k-r}\end{array}\right).\lb{m5}\eea
By Remark 3.1 we have
   \bea m^*(A^TC)=m^*\left( \left(\begin{array}{cc}0& A_1^TB_3\\ B_3^TA_1& B_1^TB_3+ D_1^TD_3\end{array}
       \right)\right),\qquad *=+,-,0.\lb{m6}\eea

Since $0\le\rank B_3=\lm\le\min\{r,\,k-r\}$, there exist $r\times r$
and $(k-r)\times (k-r)$ real invertible matrices $G_1$ and $G_2$
such that
  \bea
  G_1B_3G_2=\left(\begin{array}{cc}I_\lm&0\\0&0\end{array}\right):=F.\lb{m0}\eea
Note that if $\lm=0$ then $B_3=0$, if $\lm=\min\{r,\,k-r\}$ then
$G_1B_3G_2=\left(\begin{array}{cc}I_\lm&0\end{array}\right)\;{\rm
or}\; \left(\begin{array}{c}I_\lm\\0\end{array}\right)$, if
$\lm=r=k-r$ then $G_1B_3G_2=I_\lm$. The proof below can still go
through by corresponding adjustment.

\noindent By (\ref{m0}) we have
  \bea &&\left(\begin{array}{cc}G_1 A_1^{-1}&0\\0&G_2^T\end{array}\right)
         \left(\begin{array}{cc}0& A_1^TB_3\\ B_3^TA_1& B_1^TB_3+ D_1^TD_3\end{array}\right)
         \left(\begin{array}{cc}A_1^{-1}G_1^T&0\\0&G_2\end{array}\right)\nn\\
       &=&\left(\begin{array}{cc}0& G_1B_3G_2\\ G_2^TB_3^TG_1^T& U\end{array}\right)=\left(\begin{array}{cccc}0&0& I_\lm&0\\0&0&0&0\\
         I_\lm&0& U_1& U_2 \\0&0& U_2^T& U_4\end{array}\right).\lb{m7}\eea
Then
    \bea &&\left(\begin{array}{cccc}I_\lm&0&0&0\\0&I_{r-\lm}&0&0\\
        -\frac{1}{2}U_1&0&I_\lm&0\\-U_2^T&0&0&I_{k-r-\lm}\end{array}\right)
        \left(\begin{array}{cccc}0&0& I_\lm&0\\0&0&0&0\\
         I_\lm&0& U_1& U_2 \\0&0& U_2^T& U_4\end{array}\right)
        \left(\begin{array}{cccc}I_\lm&0&-\frac{1}{2}U_1&-U_2\\0&I_{r-\lm}&0&0\\
        0&0&I_\lm&0\\0&0&0&I_{k-r-\lm}\end{array}\right)\nn\\
      &=& \left(\begin{array}{cccc}0&0& I_\lm&0\\0&0&0&0\\
        I_\lm&0&0&0\\0&0&0& U_4\end{array}\right).\lb{m8}\eea
 Set \bea q^*=m^*(U_4),\quad *=\pm,\,0\lb{q3}\eea Then
$q^++q^0+q^-=k-r-\lm$ and (\ref{m2})-(\ref{m4}) hold from
(\ref{m6}), (\ref{m7}) and (\ref{m8}).

Also by (\ref{m8}) and Lemma 3.1 we have
 \bea \left(\begin{array}{cc}I_{k-r-\lm}&0\\2U_4&
 I_{k-r-\lm}\end{array}\right)\approx
 \left(\begin{array}{cc}1&1\\0&1\end{array}\right)^{\diamond
 q^-}\diamond I_2^{\diamond q^0}\diamond \left(\begin{array}{cc}1&-1\\0&1\end{array}\right)^{\diamond
 q^+}.\lb{q4}\eea
 By (\ref{m7}) we have
 \bea
 &&\diag((G^T_1)^{-1}A_1,G_2^{-1},G_1A_1^{-1},G_2^T)(N_kR_1^{-1}N_kR_1)\diag(A_1^{-1}G^T_1,G_2,A_1G_1^{-1},(G_2^T)^{-1})\nn\\
 &=&\left(\begin{array}{cccc}I_r&2E&2\td{A}_1&0\\0&I_{k-r}&0&0\\0&2F&I_r&0\\2F^T&2U&2E^T&I_{k-r}\end{array}\right):=M,
 \lb{m11}\eea
where $\td{A}_1=(G_1^T)^{-1}A_1  G_1^{-1} $,
$E=(G_1^T)^{-1}A_1B_3G_2=\td{A}_1F$.

Since $M$ is symplectic, we have $M^TJ_kM=J_k$. Then we have
$E=\td{A}_1F$. Since $\td{A}_1=(G_1^T)^{-1}A_1  G_1^{-1} $, it is
congruent to $\diag(a_1,a_2,\dots,a_r)$  with
 \bea a_i=1,\;1\le i\le p
\;\quad{\rm and}\;\quad  a_j=-1,\;p+1\le j\le r \;  {\rm for\;
some} \; 0\le p\le r. \lb{q5}\eea
 Then there is an
invertible $r\times r$ real matrix $Q$ such that $\det Q>0$ and
 \bea Q\td{A}_1Q^T=\diag(a_1,a_2,\dots,a_r)=\diag(\diag(a_1,a_2,...,a_\lm),\diag(a_{\lm+1},...,a_r)):=
 \diag(\Lm_1,\Lm_2).\lb{pp2}\eea
Since $\det Q>0$ we can joint it to $I_r$ by invertible continuous
matrix path. So there is a continuous invertible symmetric matrix
path $\beta$ such that
 $\alpha_1(1)=\td{A}_1$ and $\alpha_1(0)=\diag(a_1,a_2,\dots,a_r)$ with
\bea m^*(\alpha_1(t))=m^*(\td{A}_1)=m^*(A_1),\quad t\in[0,1],\;
 *=+,-.\nn\eea
Define symmetric matrix path
  \bea
  \alpha_2(t)=\left(\begin{array}{cc}2tU_1&2tU_2\\2tU^T_2&2U_4\end{array}\right),\quad
  t\in[0,1].\nn\eea

  For $t\in[0,1]$, define
 \bea \beta(t)=\left(\begin{array}{cccc}I_r&2\alpha_1(t)F&2\alpha_1(t)&0\\0&I_{k-r}&0&0
 \\0&2F&I_r&0\\2F^T&\alpha_2(t)&2F^T\alpha_1(t)^T&I_{k-r}\end{array}\right).\nn\eea
Then since $M$ is symplectic, it is easy to check that $\beta$ is a
continuous symplectic matrix path. Since
$F=\left(\begin{array}{cc}I_\lm&0\\0&0\end{array}\right)$, and
$\alpha_1(t)$ is invertible, by direct computation,
     we have
    \bea
    \rank(\beta(t)-I_{2k})=2\lm+\rank(\alpha_1(t))+\rank(U_4)\nn
                          =2\lm+r+m^+(U_4)+m^-(U_4).\nn\eea
 Hence \bea \nu_1(\beta(t))=\nu_1(\beta(1))=\nu_1(M),\qquad
t\in[0,1].\nn\eea} Since $\sg(\beta(t))=\{1\}$,   by Definition
2.2 and Lemma 2.1 \bea M&=&\beta(1)\approx
 \beta(0)\nn\\
  &=&  \left(\begin{array}{cccccc}I_\lm&0&2\Lm_1&2\Lm_1&0&0\\
      0&I_{r-\lm}&0&0&2\Lm_2&0\\
      0&0&I_\lm&0&0&0\\0&0&2I_\lm&I_\lm&0&0\\
      0&0&0&0&I_{r-\lm}&0\\
      2I_\lm&0&0&2\Lm_1&0&I_\lm\end{array}\right)\diamond\left(\begin{array}{cc}
       I_{k-r-\lm}&0\\2U_4&I_{k-r-\lm}\end{array}\right)\nn\\
 &\approx & \left(\begin{array}{cccc}I_\lm& 2\Lm_1 &2\Lm_1&0\\
      0&I_\lm&0&0\\0&2I_\lm&I_\lm&0\\2I_\lm&0&2\Lm_1&I_\lm\end{array}\right)\diamond\left(\begin{array}{cc}
       I_{r-\lm}&2\Lm_2\\0&I_{r-\lm}\end{array}\right)\diamond\left(\begin{array}{cc}
       I_{k-r-\lm}&0\\2U_4&I_{k-r-\lm}\end{array}\right)\nn\\
 &=&\left(\begin{array}{cccc}I_\lm&2\Lm_1 &2\Lm_1&0\\
      0&I_\lm&0&0\\0&2I_\lm&I_\lm&0\\2I_\lm&0&2\Lm_1&I_\lm\end{array}\right)\diamond
       \Diamond_{j=\lm+1}^r\left(\begin{array}{cc}1&2a_j\\0&1\end{array}\right)\diamond
      \left(\begin{array}{cc}I_{k-r-\lm}&0\\2U_4&I_{k-r-\lm}\end{array}\right).\lb{m9}\eea
We define continuous symplectic matrix path
     \bea \psi(t)=\left(\begin{array}{cccc} I_\lm& 2(1- t^2 )\Lm_1&2\Lm_1&0\\
      0&(1+t)I_\lm&0&0\\0&2(1- t^2)I_\lm&I_\lm&0\\2(1-t)I_\lm&0&2(1-t)\Lm_1&\frac{1}{1+t}I_\lm\end{array}\right), \quad t\in[0,1].\nn\eea
Since $\Lm_1$ is invertible, we have $\nu (\psi(t))\equiv \lm$ for
$t\in[0,1]$. So by $\sg(\psi(t))\cap \U=\{1\}$ for $t\in[0,t]$ and
Definition 2.2 we have
     \bea \left(\begin{array}{cccc}I_\lm&\Lm_1&2\Lm_1&0\\
      0&I_\lm&0&0\\0&2I_\lm&I_\lm&0\\2I_\lm&0&2\Lm_1&I_\lm\end{array}\right)=\psi(0)&\approx&\psi(1)
      =\left(\begin{array}{cc}I_\lm&2\Lm_1\\0&I_\lm\end{array}\right)\diamond
      \left(\begin{array}{cc}2I_\lm&0\\0&\frac{1}{2}I_\lm\end{array}\right)\nn\\
    &=&\Diamond_{j=1}^\lm\left(\begin{array}{cc}1&2a_j\\0&1\end{array}\right)\diamond
      D(2)^{\diamond \lm}.\lb{m10}\eea
Then by (\ref{m9}), (\ref{m10}) and Remark 2.1 we have
 \bea M \approx \left(\Diamond_{j=1}^r\left(\begin{array}{cc}1&a_j\\0&1\end{array}\right)\right) \diamond
      D(2)^{\diamond \lm}\diamond\left(\begin{array}{cc}I_{k-r-\lm}&0\\
      U_4&I_{k-r-\lm}\end{array}\right).\lb{m14}\eea
  So by (\ref{q4}), (\ref{q5}) and Remark 2.1, we have
  \bea  M\approx  \left(\begin{array}{cc}1&1\\0&1\end{array}\right)^{\diamond(p+q^-)}
  \diamond \left(\begin{array}{cc}1&-1\\0&1\end{array}\right)^{\diamond (r-p+q^+)
  }\diamond I_2^{\diamond q^0}\diamond
      D(2)^{\diamond \lm}.\lb{q6}\eea

 By Lemma 3.2, (\ref{m12}) and (\ref{m11}), we have
     \bea N_kR^{-1}N_kR\approx M.\lb{m13}\eea
Then (\ref{m1}) holds from   (\ref{q6}) and (\ref{m13}). The
proof of Lemma 3.8 is complete.\hfill\hb

 \setcounter{equation}{0}
\section {The mixed $(L_0,L_1)$-concavity}

 \noindent{\bf Definition 4.1.} {\it The mixed
$(L_0,L_1)$-concavity and mixed $(L_1,L_0)$-concavity of a
symplectic path $\gamma\in \mathcal{P}_{\tau}(2n)$ are defined
respectively by}
$$\mu_{(L_0,L_1)}(\ga)=i_{L_0}(\ga)-\nu_{L_1}(\ga),\;\;\mu_{(L_1,L_0)}(\ga)=i_{L_1}(\ga)-\nu_{L_0}(\ga).$$

By by Proposition C of \cite{LZZ}, Proposition 6.1 of
\cite{LiuZhang} and Theorem 3.1, we have the following result.

\noindent{\bf Proposition 4.1.} {\it There hold}
\bea \mu_{(L_0,L_1)}(\ga)+\mu_{(L_1,L_0)}(\ga)&=& i(\ga^2)-\nu(\ga^2)-n,\lb{pp3}\\
 \mu_{(L_0,L_1)}(\ga)-\mu_{(L_1,L_0)}(\ga)&=&concav^*_{(L_0,L_1)}(\ga)=\frac 12 \sgn M_{\varepsilon}(\ga (\tau)), \;\;0<-\varepsilon\ll 1.
\lb{pp4}\eea

Theorem 1.3 in section 1 is a special case of the following result.

\noindent{\bf Theorem 4.1.} {\it For $\ga\in \P_\tau(2n)$, let
$P=\ga(\tau)$. If $i_{L_0}(\gamma)\ge 0$, $i_{L_1}(\gamma)\ge 0$, $i(\ga)\ge n$,
$\ga^2(t)=\ga(t-\tau)\ga(\tau)$ for all $t\in[\tau,2\tau]$,  then
 \bea
 \mu_{(L_0,L_1)}(\ga)+S_{P^2}^+(1)\ge 0, \lb{dd0} \\
 \mu_{(L_1,L_0)}(\ga)+S_{P^2}^+(1)\ge 0. \lb{dd1}\eea}

{\bf Proof.} The proofs of (\ref{dd0}) and (\ref{dd1}) are almost
the same. We only prove (\ref{dd1}) which  yields Theorem 1.3.

\noindent{\bf Claim 4.1.} Under the conditions of Theorem 4.1, if
 \bea
P^2\approx
\left(\begin{array}{cc}1&1\\0&1\end{array}\right)^{\diamond
p_1}\diamond D(2)^{\diamond p_2}\diamond \td{P},\lb{m19}\eea then
    \bea i(\ga^2)+2S_{P^2}^+(1)-\nu(\ga^2)\ge n+p_1+p_2.\lb{m16}\eea
{\bf Proof of Claim 4.1.} By Theorem 7.8 of \cite{Long0} we have
     \bea P &\approx& I_2^{\diamond q_1}\diamond
     \left(\begin{array}{cc}1&1\\0&1\end{array}\right)^{\diamond
     q_2}\diamond \left(\begin{array}{cc}1&-1\\0&1\end{array}\right)^{\diamond
     q_3}\diamond (-I_2)^{\diamond q_4}\diamond \left(\begin{array}{cc}-1&1\\0&-1\end{array}\right)^{\diamond
     q_5}\diamond \left(\begin{array}{cc}-1&-1\\0&-1\end{array}\right)^{\diamond
     q_6}\nn\\&&\diamond R(\theta_1)\diamond\cdots \diamond R(\theta_{q_7})\diamond\cdots\diamond R(\theta_{q_7+q_8})
     \diamond N_2(\om_1,b_1)\diamond\cdots\diamond N_2(\om_{q_9},b_{q_9})\nn\\
     &&\diamond D(2)^{\diamond q_{10}}\diamond D(-2)^{\diamond q_{11}},\lb{m17}\eea
 where $q_i\ge 0$ for $1\le i\le 11$ with $q_1+q_2+\cdots+q_8+2q_9+q_{10}+q_{11}=n$, $\theta_j\in(0,\pi)$ for $1\le j\le q_7$,
 $\theta_ j\in (\pi,2\pi)$ for $q_7+1\le j\le q_7+q_8$,
 $\om_j\in(\U\setminus \R)$ for $1\le j\le q_{9}$ and
 $b_j=\left(\begin{array}{cc}b_{j1}&b{j_2}\\b_{j3}&b_{j4}\end{array}\right)$
 satisfying $b_{j2}\neq b_{j3}$ for $1\le j\le q_9$.

By (\ref{m17}) and Remark 2.1 we have
 \bea P^2&\approx& I_2^{\diamond (q_1+q_4)}\diamond
     \left(\begin{array}{cc}1&1\\0&1\end{array}\right)^{\diamond
     (q_2+q_6)}\diamond \left(\begin{array}{cc}1&-1\\0&1\end{array}\right)^{\diamond
     (q_3+q_5)}\diamond R(2\theta_1)\diamond\cdots \diamond R(2\theta_{q_7})\diamond\nn\\&&\diamond\cdots\diamond R(2\theta_{q_7+q_8})
     \diamond N_2(\om_1,b_1)^2\diamond\cdots\diamond N_2(\om_{q_9},b_{q_9})^2\diamond D(2)^{\diamond (q_{10}+q_{11})}.\lb{m18}\eea
By Theorem 7.8 of \cite{Long0} and (\ref{m19}) and (\ref{m18}) we
have
   \bea q_2+q_6\ge p_1,\qquad q_{10}+q_{11}\ge p_2.\lb{m20}\eea

 Since $\ga^2(t)=\ga(t-\tau)\ga(\tau)$ for all $t\in[\tau,2\tau]$
we have $\ga^2$ is also the twice iteration of $\ga$ in the periodic
boundary value case, so by the Bott-type formula (cf. Theorem 9.2.1
of \cite{Long1}), the proof of Lemma 4.1 of \cite{LZZ},  and
Lemma 2.2 we have
 \bea   && i(\ga^2)+2S_{P^2}^+(1)-\nu(\ga^2)\nn\\
       &=& 2i(\ga)+2S_{P}^+(1)+\sum_{\theta\in
       (0,\pi)}(S_{P}^+(e^{\sqrt{-1}\theta})\nn\\
       &&-(\sum_{\theta\in
       (0,\pi)}(S_{P}^-(e^{\sqrt{-1}\theta})+(\nu(P)- S_{P}^+(1))+(\nu_{-1}(P)-S_{P}^-(-1)))\nn\\
       &=&2i(\ga)+2(q_1+q_2)+(q_8-q_7)-(q_1+q_3+q_4+q_5)\nn\\
        &\ge& 2n+q_1+2q_2+(q_8-q_7)-(q_3+q_4+q_5)\nn\\
        &=&n+(2q_1+3q_2+q_6+2q_8+2q_9+q_{10}+q_{11})\nn\\
         &\ge& n+2q_2+q_6+q_{10}+q_{11}\nn\\
         &\ge& n+p_1+p_2, \lb{m21} \eea
where in the first equality we have used
$S^+_{P^2}(1)=S^+_P(1)+S^+_P(-1)$ and
$\nu(\gamma^2)=\nu(\gamma)+\nu_{-1}(\gamma)$, in the first
inequality we have used the condition $i(\ga)\ge n$, in the third
equality we have used that
$q_1+q_2+\cdots+q_8+2q_9+q_{10}+q_{11}=n$, in the last inequality we
have used (\ref{m20}). By (\ref{m21}) Claim 4.1 holds. \hfill\hb

Now we continue to prove Theorem 4.1. We set $\mathcal
{A}=\mu_{(L_1,L_0)}(\ga)+ S^+_{P^2}(1)(\ga)$ and
$\mathcal{B}=\mu_{(L_0,L_1)}(\ga)+ S^+_{P^2}(1)$.

By by Proposition C of \cite{LZZ} and Proposition 6.1 of
\cite{LiuZhang} we have
 \be i_{L_0}(\ga)+i_{L_1}(\ga)=i(\ga^2)-n,\quad \nu_{L_0}(\ga)+\nu_{L_1}(\ga)=\nu(\ga^2).\lb{m22}\ee
 From (\ref{m22}) or (\ref{pp3}) we have
   \bea \mathcal{A}+\mathcal{B}=i(\ga^2)+2S_{P^2}^+(1)-\nu(\ga^2)-n.\lb{m23}\eea

\noindent{\bf Case 1.} $\nu_{L_0}(\ga)=0$.

In this case  we have
    \bea i_{L_1}(\ga)+ S_{P^2}^+(1)-\nu_{L_0}(\ga)\ge 0+0-0=0.\nn\eea
  Then (\ref{dd1}) holds.

\noindent{\bf Case 2.} $\nu_{L_0}(\ga)=n$.

 In this case
 $P=\left(\begin{array}{cc}A&0\\C&D\end{array}\right)$, so $A$ is invertible and we have
  \bea m^0(A^TC)=\nu_{L_1}(P)=\nu_{L_1}(\ga).  \lb{duan1}\eea

   By Lemma 3.1 we have
 \bea
 NP^{-1}NP=\left(\begin{array}{cc}I_n&0\\2A^TC&I_n\end{array}\right)\approx
  I_2^{\diamond m^0(A^TC)}\diamond N_1(1,1)^{\diamond
  m^-(A^TC)}\diamond N_1(1,-1)^{\diamond m^+(A^TC)}.\lb{duan2}\eea
By Claim 4.1, (\ref{duan2}) and (\ref{m23}) we have
    \bea \mathcal{A}+\mathcal{B}\ge m^-(A^TC).\lb{duan3}\eea
By Theorem 3.1, Lemma 3.3 and (\ref{duan1}) we have
      \bea \mathcal{B}-\mathcal{A}\le
      n-(m^+(A^TC+\nu_{L_1}(P))=n-(m^+(A^TC)+m^0(A^TC)).\nn\eea
So we have
   \bea \mathcal{A}-\mathcal{B}\ge
   m^+(A^TC)+m^0(A^TC)-n.\lb{duan5}\eea
Then by (\ref{duan3}) and (\ref{duan5}) we have
   \bea 2\mathcal{A}\ge m^-(A^TC)+(m^+(A^TC)+m^0(A^TC))-n=0\nn\eea
which yields $\mathcal{A}\ge 0$ and (\ref{dd1}) holds.

\noindent{\bf Case 3.} $1\le \nu_{L_0}(\ga)=\nu_{L_0}(P)\le n-1$.

In this case by (i) of Lemma 3.8 we have
    \bea P:=\left(\begin{array}{cc}A&B\\C&D\end{array}\right)\sim
\left(\begin{array}{cccc}A_1&B_1&I_r&0\\0&D_1&0&0\\A_3&B_3&A_2&0\\
C_3&D_3&C_2&D_2
\end{array}\right),\nn\eea where $A_1,A_2,A_3$ are $r\times r$ matrices,
$D_1,D_2,D_3$ are $(n-r)\times (n-r)$ matrices, $B_1,B_3$ are
$r\times (n-r)$ matrices, and $C_2,C_3$ are $(n-r)\times r$
matrices. We  divide Case 3 into  the following 3 subcases.

\noindent{\bf Subcase 1.} $A_3=0$.

 In this subcase let $\lm=\rank B_3$. Then $0\le \lm\le\min\{r,n-r\}$, $A_1$ is invertible,
 $A_1A_2=I_r$ and $D_1D_2^T=I_{k-r}$, so we have $A$ is invertible, furthermore there holds
 $  m^0(A^TC)=\dim \ker C=\nu_{L_1}(P)$.
Suppose $m^+(A_1)=p$, $m^-(A_1)=r-p$, then by (iv) of Lemma 3.8 we have
  \bea  &&N_kR^{-1}N_kR\approx \left(\begin{array}{cc}1&1\\0&1\end{array}\right)^{\diamond p+q^-}
  \diamond \left(\begin{array}{cc}1&-1\\0&1\end{array}\right)^{\diamond (r-p+q^+)}
  \diamond I_2^{\diamond q^0}\diamond
      D(2)^{\diamond \lm},\lb{m28}\\
      &&m^+(A^TC)=\lm+ q^+,\lb{m29}\\ &&m^0(A^TC)=r-\lm+q^0,\lb{m30}\\&& m^-(A^TC)=\lm+q^-,\lb{m31}\eea
   where $q^*\ge 0$ for   $*=+,-,0$ and  $q^++q^0+q^-=n-r-\lm$.

Then by (\ref{m28}) and Claim 4.1 we have
       \bea i(\ga^2)+2 S_{P^2}^+(1)-\nu(\ga^2)\ge
       n+p+q^-+\lm\ge n+q^-+\lm.\lb{m33}\eea
By (\ref{m33}) and (\ref{m23}) we have
  \bea \mathcal {A}+\mathcal{B}\ge q^-+\lm.\lb{m34}\eea

  By Theorem 3.1 and Lemma 3.3, and (\ref{m29})-(\ref{m31}) we
have
   \bea &&\mathcal{B}-\mathcal{A}\nn\\&\le&
   n-m^+(A^TC)-m^0(A^TC)\nn\\&=&n-(q^++\lm+r-\lm+q^0)
   \nn\\&=&n-(r+q^++q^0).\nn\eea
So we have
  \bea \mathcal{A}-\mathcal{B}\ge
  (r+q^++q^0)-n.\lb{m35}\eea
Since $q^++q^0+q^-=n-r-\lm$, by (\ref{m34}), (\ref{m35}) we have
    \bea 2\mathcal{A}&\ge& q^-+\lm+(r+q^++q^0)-n\nn\\
         &=& (q^-+q^++q^0))-(n-r-\lm)\nn\\
         &=&0\nn\eea
which yields (\ref{dd1}).

\noindent{\bf Subcase 2.} $A_3$ is invertible.

In this case by (ii) of Lemma 3.8 we have
    \bea P\sim \left(\begin{array}{cc}A_1&I_r\\A_3&A_2\end{array}\right)\diamond
     \left(\begin{array}{cc}D_1&0\\\td{D}_3&D_2\end{array}\right):=P_1\diamond P_2,\lb{m37}\eea
     where $\td{D}_3$ is a $(k-r)\times (k-r)$ matrix.
Then by (\ref{m37}) and Lemma 3.2 we have
 \bea P^2\approx (N_rP_1^{-1}N_rP_1)\diamond
 (N_{n-r}P_2^{-1}N_{n-r}P_2).\lb{m38}\eea

\noindent Let $e(N_rP_1^{-1}N_rP_1)=2m$, by Lemma 3.7 we have $0\le
m\le r$ and
      \bea \frac{1}{2}\sgn M_\var(P_1)\le r-m,\quad 0<-\var\ll
      1.\lb{m39}\eea
Also by (\ref{m38}) and (\ref{m18}), there exists
$\td{P}_1\in\Sp(2m)$ such that
    \bea N_rP_1^{-1}N_rP_1\approx D(2)^{\diamond (r-m)}\diamond
    \td{P}_1.\lb{m40}\eea
By Lemma 3.1 we have \bea
N_{n-r}P_2^{-1}N_{n-r}P_2&=&\left(\begin{array}{cc}I_{n-r}&0\\2D_1^T\td{D}_3&I_{n-r}\end{array}\right)\nn\\&\approx&
      \left(\begin{array}{cc}1&1\\0&1\end{array}\right)^{\diamond
      m^-(D_1^T\td{D}_3)}\diamond I_2^{\diamond
      m^0(D_1^T\td{D}_3)}\diamond\left(\begin{array}{cc}1&-1\\0&1\end{array}\right)^{\diamond
     m^+(D_1^T\td{D}_3)}.\lb{m41}\eea
So by Claim 4.1 and (\ref{m40}), (\ref{m41}), (\ref{m38}) and
(\ref{m23}) we have
 \bea \mathcal{A}+\mathcal{B}\ge m^-(D_1^T\td{D}_3)+r-m.\lb{m42}\eea
By Theorem 3.1 and Lemma 3.3 together with Lemma 3.7, for
$0<-\var\ll 1$ we have
     \bea \mathcal{B}-\mathcal{A}&=&\frac{1}{2}\sgn M_\var(P_1)+
      \frac{1}{2}\sgn M_\var(P_2)\nn\\
      &\le& r-m+(n-r)-m^+(D_1^T\td{D}_3)-m^0(D_1^T\td{D}_3)\nn\\
      &=&n-(m+m^0(D_1^T\td{D}_3)+m^0(D_1^T\td{D}_3)).\nn\eea
 We remind that we have used the fact $m^0(D_1^T\td D_3)=\ker (\td D_3)=\nu_{L_1}(P_2)$. So we have
    \bea \mathcal{A}-\mathcal{B}\ge
    (m+m^+(D_1^T\td{D}_3)+m^0(D_1^T\td{D}_3))-n.\lb{m43}\eea
 Note that
  \bea
  m^+(D_1^T\td{D}_3)+m^0(D_1^T\td{D}_3)+m^-(D_1^T\td{D}_3)=n-r.\lb{m44}\eea
 Then by (\ref{m42}), (\ref{m43}) and (\ref{m44}) we have
  \bea 2\mathcal{A}&\ge&
  m^-(D_1^T\td{D}_3)+r-m+(m+m^+(D_1^T\td{D}_3)+m^0(D_1^T\td{D}_3))-n\nn\\
  &=&
  m^+(D_1^T\td{D}_3)+m^0(D_1^T\td{D}_3)+m^-(D_1^T\td{D}_3)-(n-r)\nn\\
  &=&0\nn\eea
  which yields (\ref{dd1}).

\noindent{\bf Subcase 3.} $1\le \rank A_3=l\le r-1$.

   In  this case by (iii) of Lemma 3.8 we have
    \bea P\sim
            \left(\begin{array}{cc}U&I_l\\\Lm &
       V\end{array}\right)\diamond\left(\begin{array}{cccc}
     \td{A}_1&\td{B}_1&I_{r-l}&0\\
       0&D_1&0&0\\0&\td{B}_3&\td{A}_2&0\\
       \td{C}_3&\td{D}_3&\td{C}_2&\td{D}_2\end{array}\right):=P_3\diamond P_4,\lb{m45}\eea
     where $\td{A}_1,\td{A}_2$ are $(r-l)\times (r-l)$ matrices,
     $\td{B}_1,\td{B}_3$ are $(r-l)\times (n-r)$ matrices,
     $\td{C}_2, \td{C}_3$ are $(n-r)\times (r-l)$ matrices,
     $D_1,\td{D}_2,\td{D}_3$ are $(n-r)\times (n-r)$ matrices,
     $U, V, \Lm$ are $l\times l$ matrices, and $\Lm$ is invertible.

 Let $\lm=\rank \td{B}_3$ and denote
$P_4=\left(\begin{array}{cc}\td{A}&\td{B}\\\td{C} &
       \td{D}\end{array}\right),$ where $\td{A},\td{B},\td{C},\td{D}$ are $(n-l)$-order real matrices.
  Assume $m^+(\td{A}_1)=p$, $m^-(\td{A}_1)=r-l-p$, then by (iv) of
Lemma 3.8 we have
  \bea  &&N_kP_4^{-1}N_kP_4\approx \left(\begin{array}{cc}1&1\\0&1\end{array}\right)^{\diamond (p+q^-)}
  \diamond \left(\begin{array}{cc}1&-1\\0&1\end{array}\right)^{\diamond (r-l-p+q^+)}
  \diamond I_2^{\diamond q^0}\diamond
      D(2)^{\diamond \lm},\lb{m46}\\
      &&m^+(\td{A}^T\td{C})=\lm+ q^+,\lb{m47}\\ &&m^0(\td{A}^T\td{C})=r-l-\lm+q^0,\lb{m48}
      \\&& m^-(\td{A}^T\td{C})=\lm+q^-,\lb{m49}\eea
   where $q^*\ge 0$ for $*=+,-,0$ and $q^++q^0+q^-=n-r-\lm$.

\noindent Let $e(N_lP_3^{-1}N_lP_3)=2m$, by Lemma 3.7 we have $0\le
m\le l$
      \bea \frac{1}{2}\sgn M_\var(P_3)\le l-m,\quad 0<-\var\ll
      1.\lb{m50}\eea
By similar argumet as in the proof of Subcase 2, there exists
$\td{P}_3\in\Sp(2m)$ such that
    \bea N_rP_3^{-1}N_rP_3\approx D(2)^{\diamond (l-m)}\diamond
    \td{P}_3.\lb{m51}\eea
  So by Claim 4.1, (\ref{m45}), (\ref{m46}), (\ref{m51}), and
(\ref{m23}) we have
 \bea \mathcal{A}+\mathcal{B}\ge q^-+l-m+\lm.\lb{m54}\eea
By Theorem 3.1, Lemma 3.3, (\ref{m47}), (\ref{m48}) and (\ref{m50}),
for $0\le -\var\ll 1$ we have
     \bea \mathcal{B}-\mathcal{A}&=&\frac{1}{2}\sgn M_\var(P_3)+
      \frac{1}{2}\sgn M_\var(P_4)\nn\\
      &\le& \frac{1}{2}\sgn
      M_\var(P_3)+(n-l)-m^+(\td{A}^T\td{C})-m^0(\td{A}^T\td{C})\nn\\
      &\le& l-m+(n-l)-(\lm+q^+)-(r-l-\lm+q^0)\nn\\
      &=&n+(l-m)-(q^++q^0+r).\nn\eea
So we have
    \bea \mathcal{A}-\mathcal{B}\ge
    (q^++q^0+r)-n-(l-m).\lb{m55}\eea
Since $q^++q^0+q^-=n-r-\lm$, by (\ref{m54}) and (\ref{m55}) we have
  \bea 2\mathcal{A}&\ge&
  q^-+l-m+\lm+(q^++q^0+r)-n-(l-m)\nn\\
  &=&(q^++q^0+q^-)-(n-r-\lm)
  \nn\\
  &=&0\nn\eea
 which yields (\ref{dd1}). Hence (\ref{dd1}) holds in Cases 1-3 and
 the proof of Theorem 4.1 is complete.\hfill\hb

\noindent{\bf Remark 4.1.} Both the  estimates  (\ref{dd0}) and
(\ref{dd1}) in Theorem 4.1 are optimal . In fact, we can construct a
symplectic path satisfying the conditions of Theorem 4.1 such that
the equalities in  (\ref{dd0}) and (\ref{dd1}) hold.  Let $\tau=\pi$
and $\ga(t)=R(t)^{\diamond n},\; t\in
  [0,\pi]$. It is easy to see that
  $i_{L_0}(\ga)=\displaystyle\sum_{0<t<\pi}\nu_{L_0}(\ga(t))=0$ and also $i_{L_1}(\ga)=\displaystyle\sum_{0<t<\pi}\nu_{L_1}(\ga(t))=0$,
  $\nu_{L_0}(\ga)=\nu_{L_1}(\ga)=n$, $\ga^2(t)=\ga(t-\pi)\ga(\pi)$
  for $t\in[\pi,2\pi]$, $i(\ga)=n$ and
  $P=\ga(\pi)=-I_{2n}$ hence by Lemma 2.2
  $S_{P^2}^+(1)=S_{I_{2n}}^+(1)=n$. So we have
   \bea
   \mu_{(L_0,L_1)}(\ga)+S_{P^2}^+(1)=\mu_{(L_1,L_0)}(\ga)+S_{P^2}^+(1)=0-n+n=0.\nn\eea

 \setcounter{equation}{0}
\section{Proofs of Theorems 1.1 and 1.2 }
In this section we prove Theorems 1.1-1.2.

 For $\Sg\in
\mathcal{H}_b^{s,c}(2n)$, let $j_\Sg: \Sg \rightarrow[0,+\infty)$ be
the gauge function of $\Sg$ defined by
        \bea
        j_{\Sg}(0)=0,\quad {\rm and} \quad  j_\Sg(x)=\inf\{\lambda >0\mid
        \frac{x}{\lambda}\in C\}, \quad \forall x \in
        \R^{2n}\setminus\{0\},\nn
        \eea
where $C$ is the domain enclosed by $\Sg$.

 Define
        \bea H_\alpha(x)=(j_\Sg(x))^\alpha,\;\alpha>1,\quad
        H_\Sg(x)=H_2(x),\; \forall x \in
        \R^{2n}.\lb{7.2}
        \eea
Then $H_{\Sigma} \in C^2 (\R^{2n}\backslash \{0\},\R)\cap
    C^{1,1}(\R^{2n},\R)$.

We consider the following fixed energy problem \bea
\dot{x}(t) &=& JH_\Sg'(x(t)), \lb{7.4}\\ H_\Sg(x(t)) &=& 1,   \lb{7.5}\\
x(-t) &=& Nx(t),  \lb{7.6}\\ x(\tau+t) &=& x(t),\quad \forall\,
t\in\R. \lb{7.7} \eea

 Denote by
$\mathcal{J}_b(\Sg,2)\;(\mathcal{J}_b(\Sg,\alpha)$ for $\alpha=2$ in
(\ref{7.2})) the set of all solutions $(\tau,x)$ of problem
(\ref{7.4})-(\ref{7.7}) and by $\tilde{\mathcal{J}}_b(\Sg,2)$ the
set of all geometrically distinct solutions of
(\ref{7.4})-(\ref{7.7}). By Remark 1.2 of \cite{LiuZhang} or
discussion in \cite{LZZ}, elements in $\mathcal{J}_b(\Sg)$ and
$\mathcal{J}_b(\Sg,2)$ are one to one correspondent. So we have
$^\#\td{\mathcal{J}}_b(\Sg)$=$^\#\td{\mathcal{J}}_b(\Sg,2)$.

For readers' convenience in the following we list some known results
which will be used in the proof of Theorem 1.1.

In the following of this paper, we write
$(i_{L_0}(\gamma,k),\nu_{L_0}(\gamma,k))=(i_{L_0}(\gamma^k),\nu_{L_0}(\gamma^k))$
for any symplectic path $\gamma\in \mathcal {P}_{
 {\tau}}(2n)$ and $k\in \N$, where $\ga^k$ is defined by Definition 2.5. We have

\noindent{\bf Lemma 5.1.} (Theorem 1.5 and of \cite{LiuZhang} and
Theorem 4.3 of \cite{LZ}) {\it Let $\ga_j\in \mathcal {P}_{
 {\tau_j}}(2n)$ for $j=1,\cdots,q$.  Let
 $M_j=\ga^2_j(2\tau_j)=N\ga_j(\tau_j)^{-1}N\ga_j(\tau_j)$, for $j=1,\cdots,q$. Suppose
          \bea \hat{i}_{L_0}(\ga_j)>0, \quad
          j=1,\cdots,q.\nn\eea
  Then there exist infinitely many $(R, m_1, m_2,\cdots,m_q)\in \N^{q+1}$ such that

  (i) $\nu_{L_0}(\ga_j, 2m_j\pm 1)=\nu_{L_0}(\ga_j)$,

  (ii) $i_{L_0}(\ga_j, 2m_j-1)+\nu_{L_0}(\ga_j,2m_j-1)=R-(i_{L_1}(\ga_j)+n+S_{M_j}^+(1)-\nu_{L_0}(\ga_j))$,

  (iii) $i_{L_0}(\ga_j,2m_j+1)=R+i_{L_0}(\ga_j)$.

\noindent and
    (iv) $\nu(\ga_j^2, 2m_j\pm 1)=\nu(\ga_j^2)$,

  (v) $i(\ga_j^2, 2m_j-1)+\nu(\ga_j^2,
  2m_j-1)=2R-(i(\ga_j^2)+2S_{M_j}^+(1)-\nu(\ga_j^2))$,

  (vi) $i(\ga_j^2,2m_j+1)=2R+i(\ga_j^2)$,

\noindent where we have set $i(\ga_j^2, n_j)=i(\ga_j^{2n_j})$, $\nu(\ga_j^2, n_j)=\nu(\ga_j^{2n_j})$ for $n_j\in\N$.}\\

For any $(\tau,x)\in\mathcal{J}_b(\Sg,2)$, there is a symplectic
path $\gamma_x\in {\mathcal P}_{\tau}(2n)$ corresponding to it. For
$m\in\N$, we denote by $i_{L_j}(x,m)=i_{L_j}(\ga_x^m)$ and
$\nu_{L_j}(x,m)=\nu_{L_j}(\ga_x^m)$ for $j=0,1$.  Also we denote by
$i(x,m)=i(\ga_x^{2m})$ and $\nu(x,m)=\nu(\ga_x^{2m})$. We remind
that the symplectic path $\ga_x^m$ is defined in the interval
$[0,\frac{m\tau}{2}]$ and the symplectic path $\ga_x^{2m}$ is
defined in the interval $[0,m\tau]$. If $m=1$, we denote by $i(x)=
i(x,1)$ and
$\nu(x)=\nu(x,1)$. By Lemma 6.3 of \cite{LiuZhang} we have\\

\noindent {\bf Lemma 5.2.} {\it Suppose
$^\#\tilde{\mathcal{J}}_b(\Sg)<+\infty$. Then there exist an integer
$K\ge 0$ and an injection map $\phi: \N+K\mapsto
\mathcal{J}_{b}(\Sg,2)\times \N$ such that

(i) For any $k\in \N+K$, $[(\tau,x)]\in \mathcal{J}_{b}(\Sg,2)$ and
$m\in \N$ satisfying $\phi(k)=([(\tau \;,x)],m)$, there holds
              $$i_{L_0}(x,m)\le k-1\le i_{L_0}(x,m)+\nu_{L_0}(x,m)-1,$$
where $x$ has minimal period $\tau$.

(ii) For any $k_j\in \N+K$, $k_1<k_2$, $(\tau_j,x_j)\in \mathcal
{J}_b(\Sg,2)$ satisfying $\phi(k_j)=([(\tau_j \;,x_j)],m_j)$ with
$j=1,2$ and $[(\tau_1 \;,x_1)]=[(\tau_2 \;,x_2)]$, there holds
       $$m_1<m_2.$$}

\noindent{\bf Lemma 5.3.} (Lemma 7.2 of \cite{LiuZhang}) {\it Let
$\ga\in \P_\tau(2n)$ be extended to $[0,+\infty)$ by
$\ga(\tau+t)=\ga(t)\ga(\tau)$ for all $t>0$. Suppose
$\ga(\tau)=M=P^{-1}(I_2\diamond \td{M})P$ with $\td{M}\in \Sp(2n-2)$
and $i(\ga)\ge n$. Then we have
      \bea i(\ga,2)+2S_{M^2}^+(1)-\nu(\ga,2)\ge n+2.\nn\eea}

\noindent{\bf Lemma 5.4} (Lemma 7.3 of \cite{LiuZhang}) {\it For any
$(\tau,x)\in \mathcal{J}_b(\Sg,2)$ and $m\in \N$, we have
\bea i_{L_0}(x,m+1)-i_{L_0}(x,m)&\ge& 1,\nn\\
    i_{L_0}(x,m+1)+\nu_{L_0}(x,m+1)-1&\ge&
   i_{L_0}(x,m+1)>i_{L_0}(x,m)+\nu_{L_0}(x,m)-1.\nn\eea}

{\bf Proof of Theorem 1.1.} It is suffices to consider the case
 $^\#\tilde{\mathcal{J}}_b(\Sg)<+\infty$. Since $-\Sg=\Sg$, for
 $(\tau,x) \in \mathcal{J}_b(\Sg,2)$  we have
          \bea &&H_\Sg(x)=H_\Sg(-x),\nn\\
               &&H_\Sg'(x)=- H_\Sg'(-x),\nn\\
                &&H_\Sg''(x)= H_\Sg''(-x).\lb{8.43}\eea
So $(\tau,-x)\in \mathcal{J}_b(\Sg,2)$. By (\ref{8.43}) and the
definition of $\ga_x$ we have that
 \bea \ga_x=\ga_{-x}.\nn\eea
So we have
       \bea &&(i_{L_0}(x,m),\nu_{L_0}(x,m))=(i_{L_0}(-x,m),\nu_{L_0}(-x,m)),\nn\\
       &&(i_{L_1}(x,m),\nu_{L_1}(x,m))=(i_{L_1}(-x,m),\nu_{L_1}(-x,m)),\quad \forall m\in
       \N.\lb{8.45}\eea
 So we can write
  \be \td{\mathcal {J}}_b(\Sg,2)=\{[(\tau_j,x_j)]|
j=1,\cdots,p\}\cup\{[(\tau_k,x_k)],[(\tau_k,-x_k)]|k=p+1,\cdots,p+q\}.\lb{8.46}\ee
with $x_j(\R)=-x_j(\R)$ for $j=1,\cdots,p$ and $x_k(\R)\neq
-x_k(\R)$ for $k=p+1,\cdots,p+q$. Here we remind that $(\tau_j,x_j)$
has minimal period $\tau_j$ for $j=1,\cdots,p+q$ and
$x_j(\frac{\tau_j}{2}+t)=-x_j(t), \;t\in\R$ for $j=1,\cdots,p$.

 By Lemma 5.2 we have an integer $K\ge 0$ and an injection map
 $\phi: \N+K\to
\mathcal{J}_b(\Sg,2)\times \N$. By (\ref{8.45}), $(\tau_k,x_k)$ and
$(\tau_k,-x_k)$ have the same $(i_{L_0},\nu_{L_0})$-indices.
 So by Lemma 5.2,
 without loss of generality, we can further require that
       \bea {\rm Im} (\phi)\subseteq \{[(\tau_k,x_k)]|k=1,2,\cdots,p+q\}\times
       \N.\lb{8.47}\eea
By the strict convexity of $H_\Sg$ and (6.19) of \cite{LiuZhang}),
we have
        \bea \hat{i}_{L_0}(x_k)>0,\quad
        k=1,2,\cdots,p+q.\nn\eea
Applying  Lemma 5.1 to the following associated symplectic paths
$$\ga_1,\;\cdots,\;\ga_{p+q},\; \ga_{p+q+1},\;\cdots,\;\ga_{p+2q}$$
of
$(\tau_1,x_1),\;\cdots,\;(\tau_{p+q},x_{p+q}),\;(2\tau_{p+1},x_{p+1}^2),\;\cdots,\;
     (2\tau_{p+q},x_{p+q}^2)$ respectively,
there exists a vector $(R,m_1,\cdots,m_{p+2q})\in \N^{p+2q+1}$ such
that $R>K+n$ and
   \bea &&i_{L_0}(x_k, 2m_k+1)=R+i_{L_0}(x_k),\lb{8.49}\\
        && i_{L_0}(x_k,2m_k-1)+\nu_{L_0}(x_k,2m_k-1)\nn\\
        &=&R-(i_{L_1}(x_k)+n+S_{M_k}^+(1)-\nu_{L_0}(x_k)),\lb{8.50}\eea
    for $k=1,\cdots,p+q,$ $M_k=\ga_k^2(\tau_k)$, and
     \bea &&i_{L_0}(x_k, 4m_k+2)=R+i_{L_0}(x_k,2),\lb{8.51}\\
         &&i_{L_0}(x_k,4m_k-2)+\nu_{L_0}(x_k,4m_k-2)\nn\\
         &=&R-(i_{L_1}(x_k,2)+n+S_{M_k}^+(1)-\nu_{L_0}(x_k,2)),\lb{8.52}\eea
      for $k=p+q+1,\cdots,p+2q$ and $M_k=\ga_k^4(2\tau_k)=\ga_k^2(\tau_k)^2$.

By Lemma 5.1, we also have \bea i(x_k,
2m_k+1)&=&2R+i(x_k),\lb{8.53}\\
         i(x_k,2m_k-1)+\nu(x_k,2m_k-1)
        &=&2R-(i(x_k)+2S_{M_k}^+(1)-\nu(x_k)),\lb{8.54}\eea
    for $k=1,\cdots,p+q,$ $M_k=\ga_k^2(\tau_k)$, and
     \bea i(x_k, 4m_k+2)&=&2R+i(x_k,2),\lb{8.55}\\
        i(x_k,4m_k-2)+\nu(x_k,4m_k-2)
         &=&2R-(i(x_k,2)+2S_{M_k}^+(1)-\nu(x_k,2)),\lb{8.56}\eea
      for $k=p+q+1,\cdots,p+2q$ and $M_k=\ga_k^4(2\tau_k)=\ga_k^2(\tau_k)^2$.

From (\ref{8.47}), we can set
 \bea \phi(R-(s-1))=([(\tau_{k(s)}, x_{k(s)})],m(s)),\qquad
    \forall s\in S:=\{1,2,\cdots,n\},\nn\eea
where $k(s)\in \{1,2,\cdots,p+q\}$ and $m(s)\in \N$.

We continue our proof  to study the symmetric and asymmetric orbits
separately. Let \bea S_1=\{s\in S|k(s)\le p\},\qquad S_2=S\setminus
S_1.\nn\eea
 We shall prove that
$^\#S_1\le p$ and $^\#S_2\le 2q$, together with the definitions of
$S_1$ and $S_2$, these yield Theorem 1.1.

\noindent{\bf Claim 5.1.} $^\#S_1\le p$.

\noindent {\it Proof of Claim 5.1.} By the definition of $S_1$,
$([(\tau_{k(s)},
 x_{k(s)})],m(s))$ is symmetric when $k(s)\le p$. We further prove
 that $m(s)=2m_{k(s)}$ for $s\in S_1$.

  In fact, by the definition of $\phi$ and Lemma 5.2,
   for all $s=1,2,\cdots,n$ we have
         \bea  i_{L_0}(x_{k(s)},m(s))&\le & (R-(s-1))-1=R-s \nn\\
         &\le &
         i_{L_0}(x_{k(s)},m(s))+\nu_{L_0}(x_{k(s)},m(s))-1.\lb{8.59}\eea
 By the strict convexity of $H_\Sg$ and  Lemma 2.3, we have $i_{L_0}(x_{k(s)})\ge 0$, so there holds
   \bea i_{L_0}(x_{k(s)},m(s))\le R-s< R\le R+i_{L_0}(x_{k(s)})=i_{L_0}(x_{k(s)},2m_{k(s)}+1),\lb{8.60}\eea
 for every $s=1,2,\cdots,n$, where we have used
 (\ref{8.49}) in the last equality. Note that the proofs of (\ref{8.59}) and
 (\ref{8.60}) do not depend on the condition $s\in S_1$.

It is easy to see that $\ga_{x_k}$ satisfies conditions of Theorem 4.1 with $\tau=\frac{\tau_k}{2}$. Note that by definition
$i_{L_1}(x_k)=i_{L_1}(\ga_{x_k})$ and
$\nu_{L_0}(x_k)=\nu_{L_0}(\ga_{x_k})$. So by Theorem 4.1 we have
   \be
i_{L_1}(x_k)+S_{M_k}^+(1)-\nu_{L_0}(x_k)\ge 0,\quad \forall
k=1,\cdots,p.\lb{8.64}\ee
 Hence by (\ref{8.59}) and (\ref{8.64}), if $k(s)\le p$ we have
 \bea
&&i_{L_0}(x_{k(s)},2m_{k(s)}-1)+\nu_{L_0}(x_{k(s)},2m_{k(s)}-1)-1\nn\\
&=&
R-(i_{L_1}(x_{k(s)})+n+S_{M_{k(s)}}^+(1)-\nu_{L_0}(x_{k(s)}))-1\nn\\
&\le&R-\frac{1-n}{2}-1-n\nn\\&<& R-s\nn\\
&\le&
i_{L_0}(x_{k(s)},m(s))+\nu_{L_0}(x_{k(s)},m(s))-1.\lb{8.66}\eea
 Thus
by (\ref{8.60}) and (\ref{8.66}) and Lemma 5.4 we have
 \bea 2m_{k(s)}-1< m(s)<2m_{k(s)}+1.\nn\eea
 Hence
 \bea m(s)=2m_{k(s)}.\nn\eea
So we have
 \bea \phi(R-s+1)=([(\tau_{k(s)},x_{k(s)})],2m_{k(s)}),\qquad \forall
 s\in S_1.\nn\eea
Then by the injectivity of $\phi$, it induces another injection map
 \bea \phi_1:S_1\rightarrow \{1,\cdots,p\}, \;s\mapsto
 k(s).\nn\eea
 There for $^\#S_1\le p$. Claim 5.1 is proved.

 \noindent{\bf Claim 5.2.} $^\#S_2\le 2q$.

\noindent{\it Proof of Claim 5.2.} By the  formulas
(\ref{8.53})-(\ref{8.56}), and (59) of \cite{LLZ} (also Claim 4 on
p. 352 of \cite{Long1}), we have \be m_k=2m_{k+q}\quad {\rm for}\;\;
k=p+1,p+2,\cdots,p+q.\lb{8.71}\ee
 By Theorem 4.1 we have
 \be i_{L_1}(x_k,2)+S_{M_k}^+(1)-\nu_{L_0}(x_k,2)\ge 0,\quad
p+1\le k\le p+q.\lb{8.76}\ee
 By
(\ref{8.52}), (\ref{8.59}), (\ref{8.71}) and (\ref{8.76}), for
$p+1\le k(s)\le p+q$ we have
\bea &&i_{L_0}(x_{k(s)},2m_{k(s)}-2)+\nu_{L_0}(x_{k(s)},2m_{k(s)}-2)-1\nn\\
     &=&i_{L_0}(x_{k(s)},4m_{k(s)+q}-2)+\nu_{L_0}(x_{k(s)},4m_{k(s)+q}-2)-1\nn\\
     &=&R-(i_{L_1}(x_{k(s)},2)+n+S_{M_{k(s)}}^+(1)-\nu_{L_0}(x_{k(s)},2))-1\nn\\
     &=&R-(i_{L_1}(x_k,2)+S_{M_k}^+(1)-\nu_{L_0}(x_k,2))-1-n\nn\\
     &\le&R-1-n\nn\\
     &<& R-s\nn\\
     &\le&
     i_{L_0}(x_{k(s)},m(s))+\nu_{L_0}(x_{k(s)},m(s))-1.\lb{8.77}\eea
Thus by (\ref{8.60}), (\ref{8.77}) and Lemma 5.4, we have
     \bea 2m_{k(s)}-2<m(s)<2m_{k(s)}+1,\qquad p<k(s)\le
     p+q.\nn\eea
So \bea m(s)\in \{2m_{k(s)}-1,2m_{k(s)}\}, \qquad {\rm
for}\;\;p<k(s)\le p+q.\}\nn\eea
 Especially this yields that for any $s_0$ and $s\in
S_2$, if $k(s)=k(s_0)$, then
 \bea m(s)\in
\{2m_{k(s)}-1,2m_{k(s)}\}=\{2m_{k(s_0)}-1,2m_{k(s_0)}\}.\nn\eea Thus
by the injectivity of the map $\phi$ from Lemma 5.2, we have \bea
^\#\{s\in S_2|k(s)=k(s_0)\}\le 2\nn\eea which yields Claim 5.2.

     By Claim 5.1 and Claim 5.2, we have
     \bea
     ^\#\td{\mathcal{J}}_b(\Sg)=^\#\td{\mathcal{J}}_b(\Sg,2)=p+2q\ge
     ^\#S_1+^\#S_2
     =n.\nn\eea
The proof of Theorem 1.1 is complete. \hfill\hb

{\bf Proof of Theorem 1.2.} We call a closed characteristic $x$ on
$\Sg$ a {\it dual brake orbit} on $\Sg$ if $x(-t)=-Nx(t)$. Then by
the similar proof of Lemma 3.1 of \cite{Zhang2}, a closed
characteristic $x$ on $\Sg$ can became a dual brake orbit after
suitable time translation if and only if $x(\R)=-Nx(\R)$. So by
Lemma 3.1 of \cite{Zhang2} again, if a closed characteristic $x$ on
$\Sg$ can both became brake orbits and dual brake orbits after
suitable translation, then $x(\R)=Nx(\R)=-Nx(\R)$, Thus
$x(\R)=-x(\R)$.

Since we also have $-N\Sg=\Sg$, $(-N)^2=I_{2n}$ and $(-N)J=-J(-N)$,
dually by the same proof of Theorem 1.1(with the estimate (\ref{dd0}) in Theorem 4.1), there are at least $n$
geometrically distinct dual brake orbits on $\Sg$.

If there are exactly $n$ closed characteristics on $\Sg$. By Theorem
1.1 all of them must be brake orbits on $\Sg$ after suitable time
translation. By the same argument all the $n$ closed characteristics
must be dual brake orbits on $\Sg$. Then by the argument in the
first paragraph of the proof of this theorem, all these $n$ closed
characteristics on $\Sg$ must be symmetric. Hence all of them must
be symmetric brake orbits after suitable time translation. The proof
of Theorem 1.2 is complete. \hfill\hb

\bibliographystyle{abbrv}

\begin{thebibliography}{99}
\bibitem {ABL1} A. Ambrosetti, V. Benci, Y. Long, A note on the existence of
multiple brake orbits. {\it Nonlinear Anal. T. M. A.}, 21 (1993)
643-649.

\bibitem{Be1} V. Benci, Closed geodesics for the Jacobi metric and periodic
solutions of prescribed energy of natural Hamiltonian systems. {\it
Ann. I. H. P. Analyse Nonl. } 1 (1984) 401-412.
\bibitem{BG} V. Benci, F. Giannoni, A new proof of the existence of a brake
orbit. In ``Advanced Topics in the Theory of Dynamical Systems".
{\it Notes Rep. Math. Sci. Eng.} 6 (1989) 37-49.
\bibitem{Bol} S. Bolotin, Libration motions of natural dynamical
systems. {\it Vestnik Moskov Univ. Ser. I. Mat. Mekh.} 6 (1978)
72-77 (in Russian).
\bibitem{BolZ} S. Bolotin, V.V. Kozlov, Librations with many degrees
of freedom. {\it J. Appl. Math. Mech.} 42 (1978) 245-250 (in
Russian).

\bibitem{CLM} S. E. Cappell, R. Lee, E. Y. Miller, On the
 Maslov-type index. {\it Comm. Pure Appl. Math.}, 47 (1994)
 121-186.
\bibitem {CoZ} C. Conley, E. Zehnder, Morse-type index theory for
flows and periodic solutions for Hamiltonian equations. {\it Comm.
Pure. Appl. Math.} 37 (1984), 207-253.
 \bibitem{GGP1} R. Giamb${\rm \grave{o}}$, F. Giannoni, P. Piccione,
Orthogonal Geodesic Chords, Brake Orbits and Homoclinic Orbits in
Riemannian Manifolds, {\it Adv. Diff. Eq.,} 10, (2005) 931每960.

\bibitem{GGP2} R. Giamb${\rm \grave{o}}$, F. Giannoni, P. Piccione,
 Existence of orthogonal geodesic chords on Riemannian manifolds
 with concave boundary and homeomorphic to the N-demensional disk.
  {\it Nonlinear Anal. T. M.
A.} 73 (2010) 290-337.

\bibitem{GGP3} R. Giamb${\rm \grave{o}}$, F. Giannoni, P. Piccione,
Potential wells with a unique brake orbit. Counterexamples to a
conjecture by H. Seifert. arXiv: 1203.5198vl [math. DS].

\bibitem {GZ1} H. Gluck, W. Ziller, Existence of periodic
solutions of conservtive systems. {\it Seminar on Minimal
Submanifolds,} Princeton University Press(1983) 65-98.
\bibitem {Gro} E. W. C. van Groesen, Analytical mini-max methods for Hamiltonian
brake orbits of prescribed energy. {\it J. Math. Anal. Appl.} 132
(1988) 1-12.

\bibitem {GL1}F. Guo, C. Liu,  Multiplicity of Lagrangian orbits on symmetric star-shaped hypersurfaces.
{\it Nonlinear Anal}. 69(4) (2008), 1425每1436.

\bibitem {Ha1} K. Hayashi,
Periodic solution of classical Hamiltonian systems. {\it Tokyo J.
Math.} 6(1983), 473-486.

\bibitem{HWZ} H. Hofer, K. Wysocki, and E. Zehnder, The dynamics on
three-dimensional strictly convex energy surfaces. {\it Ann. Math.}
(2) {\bf 148} (1998) 197-289.

\bibitem {Liu2} C. Liu,  Maslov-type index theory for symplectic paths with
Lagrangian boundary conditions. {\it Adv. Nonlinear Stud.} 7 (2007)
no. 1, 131--161.
\bibitem{Liu0}C. Liu, Asymptotically linear Hamiltonian systems with
 Lagrangian boundary conditions. {\it Pacific J. Math.} 232 (2007) no.1,
 233-255.
\bibitem{LLZ} C. Liu, Y. Long, C. Zhu, Multiplicity of closed
characteristics on symmetric convex hypersurfaces in $\R^{2n}$. {\it
Math. Ann.} 323 (2002) no. 2, 201--215.
\bibitem{LiuZhang} C. Liu and D. Zhang, Iteration theory of $L$-index and multiplicity of brake
orbits. arXiv: 0908.0021vl [math. SG].

\bibitem{Long0} Y. Long, Bott formula of the Maslov-type index theory.
{\it Pacific J. Math.} 187 (1999) 113-149.

\bibitem{Long1} Y. Long, Index Theory for Symplectic Paths with Applications.
Birkh\"auser. Basel. (2002).
\bibitem {Long4} Y. Long, Maslov-type index, degenerate critical
points, and asymptotically linear Hamiltonian systems. {\it Science
in China Ser. A} (1990) 673-682.

\bibitem{LZZ} Y. Long, D. Zhang, C. Zhu, Multiple brake orbits in bounded convex symmetric domains. {\it Advances in
Math.} 203 (2006) 568-635.

\bibitem {LZe}  Y. Long, E. Zehnder, Morse Theory for forced
oscillations of asymptotically linear Hamiltonian systems. In {\it
Stoc. Proc. Phys. and Geom.}, S. Albeverio et al. ed. World Sci.
(1990) 528-563.

\bibitem{LZ} Y. Long and C. Zhu, Closed characteristics on compact
convex hypersurfaces in $\R^{2n}$. {\it Ann. Math.,} {\bf 155}
(2002) 317-368.

 \bibitem {Moh} K. Mohnke, Holomorphic disks and the Chord Conjecture, {\it Ann. of Math}. 154 (2001), 219每222.

\bibitem{Ra1} P. H. Rabinowitz, On the existence of periodic solutions for a
class of symmetric Hamiltonian systems. {\it Nonlinear Anal. T. M.
A.} 11 (1987) 599-611.
\bibitem{Se1} H. Seifert, Periodische Bewegungen mechanischer Systeme.
{\it Math. Z.} 51 (1948) 197-216.
\bibitem {Sz} A. Szulkin, An index theory and existence of multiple brake orbits
for star-shaped Hamiltonian systems. {\it Math. Ann.} 283 (1989)
241-255.
\bibitem {V} C. Viterbo, A new obstruction to embedding Lagrangian
tori. {\it Invent. Math.} 100 (1990) 301-320.

\bibitem{We1} A. Weinstein, Normal modes for nonlinear Hamiltonian systems, {\it Inv. Math.} 20 (1973) 47每57.
\bibitem {brake1} D. Zhang and C. Liu,  Multiple brake orbits on compact convex
symmetric reversible hypersurfaces in $\R^{2n}$. arXiv: 1110.0722vl
[math. SG].
\bibitem {brake2} D. Zhang and C. Liu, Multiplicity of brake orbits on compact convex
symmetric reversible hypersurfaces in $\R^{2n}$ for $n\ge 4$, {\it
preprint}.
\bibitem{Zhang2} D. Zhang, Brake type closed characteristics on reversible
compact convex hypersurfaces in $\R^{2n}$.{\it Nonlinear Anal. T. M.
A.} 74 (2011) 3149-3158.
\bibitem{Zhang1} D. Zhang, Minimal period problems for brake orbits of
nonlinear autonomous reversible semipositive Hamiltonian systems.
arXiv: 1110.6915vl [math. SG].


  \end{thebibliography}

\end{document}